\documentclass{amsart}
\usepackage[dvips]{graphicx}
\usepackage{amscd}
\usepackage{amsmath}
\usepackage{amsxtra}
\usepackage{amsfonts}
\usepackage{amssymb}

\newtheorem{theorem}{Theorem}[section]
\newtheorem{corollary}[theorem]{Corollary}
\newtheorem{lemma}[theorem]{Lemma}
\newtheorem{proposition}[theorem]{Proposition}

\theoremstyle{definition}
\newtheorem{definition}[theorem]{Definition}
\newtheorem{remark}[theorem]{Remark}

\newtheorem{example}[theorem]{Example}
\theoremstyle{remark}

\renewcommand{\theclaim}{\textup{\theclaim}}

\newtheorem*{acknowledgements}{Acknowledgements}

\numberwithin{equation}{section}

\def\openone

{\mathchoice

{\hbox{\upshape \small1\kern-3.3pt\normalsize1}}

{\hbox{\upshape \small1\kern-3.3pt\normalsize1}}

{\hbox{\upshape \tiny1\kern-2.3pt\SMALL1}}

{\hbox{\upshape \Tiny1\kern-2pt\tiny1}}}

\makeatletter

\newbox\ipbox

\newcommand{\ip}[2]{\left\langle #1\, , \,#2\right\rangle}
\newcommand{\diracb}[1]{\left\langle #1\mathrel{\mathchoice

{\setbox\ipbox=\hbox{$\displaystyle \left\langle\mathstrut
#1\right.$}

\vrule height\ht\ipbox width0.25pt depth\dp\ipbox}

{\setbox\ipbox=\hbox{$\textstyle \left\langle\mathstrut
#1\right.$}

\vrule height\ht\ipbox width0.25pt depth\dp\ipbox}

{\setbox\ipbox=\hbox{$\scriptstyle \left\langle\mathstrut
#1\right.$}

\vrule height\ht\ipbox width0.25pt depth\dp\ipbox}

{\setbox\ipbox=\hbox{$\scriptscriptstyle \left\langle\mathstrut
#1\right.$}

\vrule height\ht\ipbox width0.25pt depth\dp\ipbox}

}\right. }

\newcommand{\dirack}[1]{\left. \mathrel{\mathchoice

{\setbox\ipbox=\hbox{$\displaystyle \left.\mathstrut
#1\right\rangle$}

\vrule height\ht\ipbox width0.25pt depth\dp\ipbox}

{\setbox\ipbox=\hbox{$\textstyle \left.\mathstrut
#1\right\rangle$}

\vrule height\ht\ipbox width0.25pt depth\dp\ipbox}

{\setbox\ipbox=\hbox{$\scriptstyle \left.\mathstrut
#1\right\rangle$}

\vrule height\ht\ipbox width0.25pt depth\dp\ipbox}

{\setbox\ipbox=\hbox{$\scriptscriptstyle \left.\mathstrut
#1\right\rangle$}

\vrule height\ht\ipbox width0.25pt depth\dp\ipbox}

} #1\right\rangle}

\newcommand{\cj}[1]{\overline{#1}}

\newcommand{\bz}{\mathbb{Z}}

\newcommand{\B}{\mathcal{B}}
\newcommand{\br}{\mathbb{R}}
\newcommand{\bc}{\mathbb{C}}
\newcommand{\bt}{\mathbb{T}}
\newcommand{\bn}{\mathbb{N}}

\def\blfootnote{\xdef\@thefnmark{}\@footnotetext}

\renewcommand{\mod}{\operatorname{mod}}

\input xy
\xyoption{all}
\usepackage{amssymb}

\newcommand{\Aut}[1]{\text{Aut}(#1)}

\def\-{^{-1}}
\def\B{\mathcal{B}}

\def\U{\mathcal{U}}

\def\Ind{{\operatorname*{Ind}}}
\def\Aut{{\operatorname*{Aut}}}
\def\Rep{{\operatorname*{Rep}}}

\begin{document}
\title[Dualitry for Baumslag-Solitar groups]{A duality approach to representations of Baumslag-Solitar groups}
\author{Dorin Ervin Dutkay}
\blfootnote{Research supported in part by a grant from the National Science Foundation DMS-0704191}
\address{[Dorin Ervin Dutkay]University of Central Florida\\
	Department of Mathematics\\
	4000 Central Florida Blvd.\\
	P.O. Box 161364\\
	Orlando, FL 32816-1364\\
U.S.A.\\} \email{ddutkay@mail.ucf.edu}
\author{Palle E.T. Jorgensen}
\address{[Palle E.T. Jorgensen]University of Iowa\\
Department of Mathematics\\
14 MacLean Hall\\
Iowa City, IA 52242-1419\\}\email{jorgen@math.uiowa.edu}
\thanks{} 
\dedicatory{
Dedicated to the memory of Professor George Mackey
}
\begin{abstract}
   We give an operator theoretic approach to the
constructions of multiresolutions as they are used in a number of basis
constructions with wavelets, and in Hilbert spaces on fractals. Our approach starts with the following version of the
classical Baumslag-Solitar relations  $u t  = t^2 u$ where $t$ is a unitary
operator in a Hilbert space $\mathcal H$ and $u$ is an isometry in $\mathcal H$. There are isometric dilations of this system into a bigger Hilbert space, relevant for wavelets. For a variety of carefully selected
dilations, the ``bigger'' Hilbert space may be $L^2(\br)$, and the dilated
operators may be the unitary operators which define a dyadic wavelet
multiresolutions of $L^2(\br)$ with the subspace $\mathcal H$ serving as the corresponding
resolution subspace; that is, the initialized resolution which is generated
by the wavelet scaling function(s). In the dilated Hilbert space, the
Baumslag-Solitar relations then take the more familiar form $u t u^{-1} = t^2$.

We offer an operator theoretic framework
including the standard construction; and we show how the representations of
certain discrete semidirect group products serve to classify the
possibilities. For this we analyze and compare several types of unitary representations of these semidirect products: the induced representations in Mackey's theory, the wavelet representations on $L^2(\br)$, the irreducible representation on the dual, the finite dimensional representations, and the the regular representation.

\end{abstract}
\maketitle \tableofcontents

\section{Introduction}\label{intro}
   It is interesting to observe that there is a fundamental duality in recent studies \cite{LPT01, Dut04, Dut05, Dut06, DuJo06,DuJo07} involving such diverse areas as Baumslag-Solitar groups, wavelets, wavelet-sets, dynamical systems, and $C^*$-algebra crossed products. While the duality uses ideas from Mackey's imprimitivity analysis and Frobenius reprocity, the contexts where we use them are more subtle than a traditional and simpler context of Mackey's methods for Lie groups. As we demonstrate in the present paper, there are several reasons for this; the presence of non-type $I$ representations is just one of them.

   A crucial notion in Mackey's approach to representations of semidirect products of continuous groups is measurable cross-section, see e.g., \cite{Mac49, Mac63, Mac76}. However as we show in the main part of our paper, things are very different for our present discrete semidirect products.  This affects both our application of this non-abelian harmonic analysis, as well as our decompositions of data (in the form of functions) into its spectral constituents. As background references to Mackey cross-sections and operator algebras, we give \cite{Mac63}, \cite{Arv76} and \cite{KaRi86}.

   Much of the current and systematic representation theory for non-type $I$ groups begins with Thoma's paper \cite{Tho64}. This further inspired more connections, in the spirit of G.W. Mackey, between groups and ergodic theory; see for example the books \cite{Pet83} and \cite{Wal82}. Aside from these general references, we will draw here on the standard facts from duality theory for locally compact abelian groups, see e.g., \cite{HeRo63}. For general facts about crossed products for groups and operator algebras, and their ideal structure, the reader may wish to consult Williams et al \cite{CMW84, Wil82}.

       In our planning of this contribution, two recent and related papers inspired us the most: They are \cite{LPT01} and \cite{MaVa00}. Especially \cite{LPT01} points out the intriguing discovery that wavelet sets (Definition \ref{def5_1}) as introduced in \cite{DaLa98} arise as support sets for direct integral decompositions of the groups under study; an observation which surely would have pleased G. W. Mackey. 

     In section \ref{plancherel}, we show that the adjoint representation for the normal abelian subgroup $B$ in $G$ is simply the action of $\alpha$ (or rather integral powers of $\alpha$), and that the co-adjoint action of $G$ on $K = \hat B$ (Pontryagin dual) is the action by the dual automorphism $\hat \alpha$. Our Proposition \ref{prop3_13} below states that this version of Mackey's  co-adjoint action in fact does {\it not} have a measurable cross-section, and we have non-type $I$ representations. For the benefit of readers, and for comparison, we outline in section \ref{concluding} the simplest case \cite{Mac49} of standard Mackey orbit theory, the one motivated by the Stone-von Neumann uniqueness theorem.

       Before turning to our results we outline the framework. Our starting point is an infinite discrete abelian group $B$, and a fixed automorphism $\alpha$ of $B$. By Pontryagin duality, the dual group $K := \hat B$ of all unitary characters $\chi$ on $B$ is compact. The group $K$ carries a dual automorphism $\hat\alpha$. In the applications mentioned above, several versions of Fourier duality will be involved, and  we will have occasion to work with two crossed products, one from the pair $(B, \alpha)$ and the other from the dual pair $(K, \hat\alpha)$. The first will be denoted $G := B \rtimes_\alpha \bz$, and the second $L := K \rtimes _{\hat\alpha}\bz$. The second will play a role in the analysis of the unitary representations of the first. Moreover the groups that arise as $G$ include the traditional Baumslag-Solitar groups.

      Since every element $\chi$ in $K$ is a one-dimensional representation of $B$, in the study of representations of $G$, it is natural to turn to the induced representations $\Ind_B^G(\chi)$. In our first result we show that these induced representations are always infinite-dimensional, but their commutants depend on whether or not $\chi$ has a finite period.

    Nonetheless, we show that the so-called scaling generator in $\Ind_B^G(\chi)$ as a unitary operator always has the same spectral type (Lebesgue spectrum), regardless of what properties the inducing character $\chi$ has.

     Moreover, we show that the induced representations $\Ind_B^G(\chi)$ is irreducible iff $\chi$ has infinite period (i.e., is ``aperiodic'').

      Even if $\chi$ has finite period, the scaling generator in $\Ind_B^G(\chi)$ as a unitary operator in fact has the same spectral type as translation of bi-infinite $l^2$-sequences, so the bilateral shift in $l^2(\bz)$; the bi-lateral shift of multiplicity one, which of course in turn is unitarily equivalent with multiplication by $z$ in the $L^2(\bt)$ picture.

     If $\chi$ has finite period $p$ say, then it is naturally associated with it a $p$-dimensional irreducible representation of $G$, but we show that this representation is not induced.  The scaling generator of this irreducible representation as a unitary operator is the cyclic permutation of $\bz_p = \bz/p\bz = \{0, 1,\dots, p - 1\}$, i.e, on the natural basis, the operator $P$ sends $0$ to $p - 1$, $1$ to $0$, $2$ to $1$, $3$ to $2$, etc, $p - 1$ to $p - 2$.

     As an operator on $l^2(\bz_p)$, $P$ is of course a unitary operator. Even though in this case the induced representation $\Ind_B^G(\chi)$ is reducible, we show that it does not include the irreducible $p$-representation.

     The disjointness of the two classes of representations is reflected in
the unitary operators that represent the scaling part in the semidirect
product group $G$. For one, this operator $T$ is the bilateral shift in $l^2(\bz)$,
and for the other it is the cyclic permutation $P$ of the basis vectors in
$l^2(\bz/p\bz)$. The deeper reason for why the irreducible $p$-representation of $G$ can not be boundedly embedded in
$\Ind_B^G(\chi)$ (even if $\chi$ has period $p$) is that $P$ can not be a matrix corner in $T$. This can be seen for example from an application of the Riemann-Lebesgue lemma.

     If the inducing character $\chi$ has a finite period $p$, then the induced representation $\Ind_B^G(\chi)$ has a direct integral decomposition.

To work out the direct integral decompositions we will give several applications of our result on its commutant. Specifically, we show that the commutant of the induced representation $\Ind_B^G(\chi)$ is the algebra of all the multiplication operators  $\{f(z^p) \,|\, f \in L^\infty(\bt)\}$ acting on the usual Hilbert space $L^2(\bt)$, $\bt = \br/\bz$.

So the projections in the commutant of $\Ind_B^G(\chi)$ come from $f =  f_S =$ an indicator function of a measurable subset $S$ of $\bt$,  $\{f(z^p)\, |\, f \in L^\infty(\bt)\}$.

    A good part of our present motivation derives from a recent and inspiring paper \cite{LPT01}. It offers a fascinating new insight into the analysis of wavelets, and it is based in a surprising and intriguing way on Mackey's theory of induced representations. While it may be thought of as an application of decomposition theory for unitary representations to wavelets, in fact, deep results on wavelets work also in the other direction: wavelet analysis throws new light on a fundamental question in representation theory.

        Our references to Mackey theory and induced representations are \cite{Mac76, Jor88, Ors79}, and to wavelet analysis \cite{Dut04, Dut05, Dut06,DuJo06, DuJo07}. In addition we will make use of the fundamentals on {\it wavelet sets}, see especially \cite{LPT01, DLS98}. In our construction we will further make use of facts from the theory of crossed products from \cite{BrJo91} and \cite{KTW85}.

\section{Wavelet sets}

The contrast between the decomposition issues in the present mix of
cases of continuous and discrete groups is illustrated nicely for what in
the wavelet literature is called wavelet sets. They have both analytic and
geometric significance, see \cite{DLS98} and \cite{BJMP05, Mer05}.

        It was recently discovered by Larson et al that there is a specific class of monic wavelets in $\br^d$ for all $d$; and that via the Fourier transform in $L^2(\br^d)$ they are parameterized by certain measurable subsets $E$ in $\br^d$ which tile $\br^d$ under two operations: one under translations by the unit rank-$d$ lattice $\bz^d$, and the other by transformations under a certain expansive $d$ by $d$ matrix $A$ over $\bz$, i.e., under the scaling transformations $A^j$ as $j$ ranges over $\bz$, i.e., under stretching and squeezing of sets in $\br^d$ under powers of a fixed matrix.

       On the other hand, there is a class of discrete semidirect product groups $G$ generated by the same pair $(\bz^d, A)$; and Mackey's theory lends itself naturally to the study of these groups $G$. In fact by Mackey induction, there is a family of monomial representations of $G$ naturally indexed by points $\chi$ in $K=\hat B$. But in general, we further know that wavelets are also derived from a certain canonical unitary representation $U_w$ of $G = G(\bz^d, A)$ acting by operators in the Hilbert space $L^2(\br^d)$, and the result in \cite{LPT01} is that there is a one-to-one correspondence between wavelet sets $E$ on one side, and sets $E$ which support an orthogonal direct integral decomposition for the representation $U_w$ on the other. Since it is known that wavelet sets may be disjoint, it follows in particular that $U_w$ may have direct integral decompositions in the sense of Mackey with support on disjoint Borel subsets.

       In particular, we show that this phenomenon occurs naturally, and in an especially constructive manner. The earlier known instances, e.g., \cite{Mac76}, of such multiplicity or dichotomy for sets that support direct integral decomposition have been rather abstract in nature, or rather this was one of the first examples of two inequivalent direct integral decompositions.

\section{Spectral types}\label{spectral}

While the initial Baumslag-Solitar operator relations have an isometric scaling operator, we look for useful unitary dilations. It turns out that there are two candidates for the corresponding spectral types: Starting with a finite rank lattice, we get an extended discrete abelian group $B$, and automorphism $\alpha$ in $B$, and a semidirect product $G = B \rtimes_\alpha \bz$. The compact dual $K=\hat B$ carries a dual automorphism $\hat \alpha$.

In this section, we establish a mapping from certain orbits in $\hat B$ into a class of induced representations, in fact mapping into equivalence classes of representations of the group $G$, with representations induced from points $\chi$ in $\hat B$ being mapped into irreducible representations of $G$.

We prove two things:

(1) The mapping is onto a class of continuous spectrum representations.

(2) We show which representations are associated with which wavelets.

Since representations $U$ of $G$ are determined by $G$-orbits in $K:=\hat B$, the spectral type of the corresponding restrictions $U|_B:b\rightarrow U(b)$ to $B\subset G$ is an important invariant. For $F\in\mathcal H$, let $\mu_F$ be the spectral measure of $U|_B$, i.e.,
$$\ip{F}{U(b)F}=\int_K\chi(b)\,d\mu_F(\chi)=\hat\mu_F(b).$$
Set $$\mathcal {H}_p=\mathcal H_p(U):=\{F\in\mathcal H\,|\, \mu_F\mbox{ is atomic}\}.$$
Pick an invariant mean $m_B$ on $B$. Then by Wiener's lemma
$$\mathcal H_p=\{F\in\mathcal H\,|\, m_B(|\ip{F}{U(b)F}|^2)>0\};$$
and
$$m_B(|\ip{F}{U(b)F}|^2)=\sum_{\mbox{atoms}}|\mu_F(\{\mbox{atoms}\})|^2.$$

But the spectrum of the $B$-restricted representations may be discrete, or not.  The absence of atoms (when the measure is continuous) is decided by a certain mean of the square expression, as described above.

The vectors in $\mathcal H$ for which the expression is $0$, or for which it is $> 0$, form closed subspaces which reduce the unitary representation U.

If $\mathcal H_p = \mathcal H$, then $U$ is induced from some $\chi$.

If $\mathcal H_p = 0$, then $U$ is disjoint from every induced representation $\Ind_B^G(\chi)$.

\begin{definition}\label{def2_1}
We now turn to definitions and basic facts. Given: 
\begin{itemize}
\item[$\bullet$] $B$: a fixed discrete abelian group;
\item[$\bullet$] $\alpha\in\Aut(B)$ a fixed automorphism of $B$;
\item[$\bullet$] $K:=\hat B=$the Pontryagin dual, i.e., the group of all unitary characters on $B$: $\chi:B\rightarrow\bt=\{z\in\bc\,|\,|z|=1\}$ such that $\chi(b+c)=\chi(b)\chi(c)$, for all $b,c\in B$;
\item[$\bullet$] $\hat\alpha\in\Aut(K)$ denotes the dual automorphism, i.e., $(\hat\alpha\chi)(b)=\chi(\alpha(b))$, $\chi\in K,b\in B$.
\end{itemize}
\end{definition}

\begin{definition}\label{def2_2}
{\it Semidirect products:} $G:=B\rtimes_\alpha\bz$ will denote the semidirect product of $B$ with the automorphism $\alpha$, i.e., 
\begin{equation}\label{eq1}
(j,b)(k,c)=(j+k,\alpha^j(c)+b),\quad(j,k\in\bz,b,c\in B).
\end{equation}

\end{definition}

\begin{example}\label{ex2_3}
The simplest example of this type of relevance to wavelet analysis is the following:
\begin{equation}\label{eq3}
B:=\bz\left[\frac12\right]:=\bigcup_{k\geq0}2^{-k}\bz\subset\br.
\end{equation}
Note that  $\bz\subset\frac12\bz\subset\dots\subset\frac1{2^k}\bz\subset\dots$ so $B$ is a subgroup of $(\br,+)$, and it is an inductive limit of the rank-one groups $2^{-k}\bz$, $k=0,1,2,\dots$.
\par
Note however that we use the discrete topology on $\bz\left[\frac12\right]$ and not the Euclidean topology induced from $\br$.
\par
A direct check reveals that
\begin{equation}\label{eq4}
\alpha(b):=2b,\quad (b\in B),
\end{equation}
defines an automorphism of $B$.
\par
It is well known that $\hat \br\cong\br$ with ``$\hat{}$'' referring to Pontryagin duality. From \eqref{eq3} we conclude that there is a natural embedding
\begin{equation}\label{eq5}
\br\hookrightarrow K,\quad t\mapsto\chi_t
\end{equation}
with dense range, often referred to as an infinitely winding curve on the ``torus'' $K$. Specifically,
\begin{equation}\label{eq6}
\chi_t(b)=e^{i2\pi tb},\quad(b\in B,t\in\br).
\end{equation}

In general, points in $K_2:=\widehat{\bz\left[\frac12\right]}$ will be written as infinite words
\begin{equation}\label{eq7}
(z_0,z_1,z_2,\dots),\quad z_k\in\bt,z_{k+1}^2=z_k, k\in\bz, k\geq 0.
\end{equation}
Then
\begin{equation}
\hat\alpha(z_0,z_1,z_2,\dots)=(z_0^2,z_0,z_1,\dots),\mbox{ and }\hat\alpha^{-1}(z_0,z_1,z_2,\dots)=(z_1,z_2,z_3,\dots).
\end{equation}
\end{example}

\begin{example}\label{ex2_4}
Our discussion of the example $B=\bz\left[\frac12\right]$ carries over {\it mutatis mutandis} to the following class of structures:
\par
Given $d\in\bn$. Let $A$ be a $d\times d $ matrix over $\bz$ with $\det A\neq0$. Let $\bz^d$ be the standard rank-$d$ lattice realized as a subgroup in $\br^d$. Then 
let $B:=\bz_d[A^{-1}]$ be the inductive limits of the groups 
\begin{equation}\label{eq9}
\bz^d\hookrightarrow A^{-1}\bz^d\hookrightarrow A^{-2}\bz^d\hookrightarrow\dots\hookrightarrow A^{-k}\bz^d\hookrightarrow\dots
\end{equation}
or equivalently
\begin{equation}\label{eq10}
\bz_d[A^{-1}]:=\bigcup_{k\geq0} A^{-k}\bz^d.
\end{equation}
Generalizing \eqref{eq4}, we get
\begin{equation}\label{eq11}
\alpha(b):=Ab,\quad (b\in B)
\end{equation}
where $Ab$ refers to matrix multiplication, i.e., 
\begin{equation}\label{eq12}
(Ab)_j=\sum_{k=1}^dA_{j,k}b_k
\end{equation}
where $b$ is  viewed as a column vector. 
\end{example}

We now return to the general case of Definitions \ref{def2_1} and \ref{def2_2}.

First note that generally, the group $G=B\rtimes_\alpha\bz$ is discrete, and that for both $G$ and the subgroup $B\leftrightarrow\{(0,b)\,|\,b\in B\}\subset G$, the respective Haar measures are simply the counting measure.

Motivated by wavelets, we are interested in the unitary representations of $G$, i.e. 
\begin{equation}\label{eq13}
U:G\rightarrow \mathcal B(\mathcal H),\mbox{ such that }U(g_1g_2)=U(g_1)U(g_2),\, U(g)^*=U(g^{-1}),\quad(g,g_1,g_2\in G),
\end{equation}
where $\mathcal H$ is some Hilbert space. Here $\mathcal B(\mathcal H)$ denotes the algebra of all bounded operators on $\mathcal H$.

Since $B$ is abelian, then, by Stone's theorem, for every unitary representation $V\in\Rep(B,\mathcal H_V)$ there is a projection valued measure $P$ defined on the Borel subsets in the dual group $K$ such that
\begin{equation}\label{eq14}
V(b)=\int_K\chi(b)\,dP(\chi),\quad(b\in B).
\end{equation}
Here $V$ is the direct integral of one-dimensional representations, i.e., points $\chi$ in $K$, or equivalently $\mathcal H_\chi=\bc$.

\begin{definition}\label{def2_5}
{\it Representations of $G$ which are induced from points in $K$.} 

Let $\chi\in K$ be given, and set $\mathcal F=\mathcal F(\chi):=$ all measurable functions $F:G\rightarrow\bc$ such that 
\begin{equation}\label{eq15}
F(j,b)=\chi(b)F(j,0),\quad(b\in B,j\in\bz)
\end{equation}
and
\begin{equation}\label{eq16}
\|F\|_\chi^2:=\sum_{j\in\bz}|F(j,0)|^2<\infty.
\end{equation}

It is imediate that $\mathcal F(\chi)$ is a Hilbert space $\mathcal H(\chi)$ relative to the norm $\|\cdot\|_\chi$ in \eqref{eq16}.

The induced representation $U\in\Rep(G,\mathcal H(\chi))$, $U:=\Ind_B^G(\chi)$ is defined by
\begin{equation}\label{eq17}
(U(g_1)F)(g_2):=F(g_2g_1),\quad(g_1,g_2\in G, F\in\mathcal H(\chi)).
\end{equation}

We leave to the reader to check that $U=\Ind_B^G(\chi)$ is indeed a unitary representation of $G$ acting in the Hilbert space $\mathcal H(\chi)$. We will be interested in representations up to unitary equivalence.
\end{definition}

\section{$\Ind_B^G(\chi)$}\label{ind}
Let $\mathcal H_2:=l^2(\bz)$ be the usual $l^2$-space of square-summable doubly-infinite sequences, i.e., $\xi=(\xi_k)_{k\in\bz}$ with norm $\|\xi\|_2^2:=\sum_{k\in\bz}|\xi_k|^2<\infty$ and inner product
\begin{equation}\label{eq18}
\ip{\xi}{\eta}:=\sum_{k\in\bz}\cj\xi_k\eta_k.
\end{equation}
If $A:\mathcal H_2\rightarrow \mathcal H_2$ is a bounded linear operator, we shall refer to its matrix $A=(A_{j,k})_{j,k\in\bz}$ as folllows
\begin{equation}\label{eq19}
A_{j,k}=\ip{\delta_j}{A\delta_k}
\end{equation}
where for all $j\in\bz$, $\delta_j(i)=\delta_{i,j}=\left\{\begin{array}{cc}
1,&\mbox{if }i=j\\
0,&\mbox{if }i\neq j\end{array}\right. .$
Note that $\{\delta_j\,|\,j\in\bz\}$ is the familiar and canonical orthonormal basis (ONB) in $\mathcal H_2$, i.e., 
$\ip{\delta_j}{\delta_k}=\delta_{j,k}$, $j,k\in\bz$, see \eqref{eq18} for the definition of $\ip{\cdot}{\cdot}$. 

If $\chi\in K=\hat B$ is given, and $(B,\alpha)$ is as in Definition \ref{def2_1}, then we set
\begin{equation}\label{eq20}
D_\chi(b)=(\chi(\alpha^k(b))_{k\in\bz}
\end{equation}
where the right-hand-side in \eqref{eq20} refers to the diagonal matrix with the specified entries, i.e., 
\begin{equation}\label{eq21}
D_\chi(b)\delta_k=\chi(\alpha^k(b))\delta_k,\quad(k\in\bz).
\end{equation}

Further, set
\begin{equation}\label{eq22}
(T_j\xi)_k:=\xi_{k+j},\quad(k,j\in\bz).
\end{equation}
It is immediate that $T$ in \eqref{eq22} defines a $T\in\Rep(\bz,l^2(\bz))$.

Under the Fourier transform
\begin{equation}\label{eq23}
f_\xi(z):=\sum_{k\in\bz}\xi_kz^k,\quad(z\in\bt),\quad l^2(\bz)\ni\xi\mapsto f_\xi\in L^2(\bt),
\end{equation}
$T$ takes the form 
\begin{equation}\label{eq24}
(T_jf)(z)=z^{-j}f(z),\quad(f\in L^2(\bt),j\in\bz).
\end{equation}

\begin{lemma}\label{lem3_1}
Let $\chi\in\hat B=:K$ be given, and for pairs $(j,b)\in G=B\rtimes_\alpha\bz$, set 
\begin{equation}\label{eq25}
U(j,b)=D_\chi(b)T_j,\quad(b\in B,j\in\bz).
\end{equation}
Then $U\in\Rep(G,\mathcal H_2)$. 
\end{lemma}

\begin{proof}
A similar fact is proved in \cite{LPT01}. We must check that
\begin{equation}\label{eq26}
U(j,b)U(k,c)=U(j+k,\alpha^j(c)+b),\quad(j,k\in\bz,b,c\in B).
\end{equation}
But this identity follows from the commutation identity
\begin{equation}\label{eq27}
T_jD_\chi(b)T_j^*=D_\chi(\alpha^j(b)),\quad(j\in\bz,b\in B);
\end{equation}
and we leave the proof of \eqref{eq27} to the reader.
\end{proof}

\begin{theorem}\label{th3_2}
Let $\chi\in\hat B=:K$ be given. Then the representation $U\in\Rep(G,\mathcal H_2)$ in \eqref{eq25} is unitarily equivalent to $\Ind_B^G(\chi)$. 
\end{theorem}

\begin{proof}
First recall the Hilbert space $\mathcal H(\chi)$ of the representation $\Ind_B^G(\chi)$ from Definition \ref{def2_5}.

For $\xi\in l^2(\bz)$, set
\begin{equation}\label{eq28}
F_\xi(j,b):=\chi(b)\xi_j,\quad(b\in B,j\in\bz).
\end{equation}
A direct verification shows that $F_\xi\in\mathcal H(\chi)$. See properties \eqref{eq15}-\eqref{eq16} in Definition \ref{def2_5}.

Setting
\begin{equation}\label{eq29}
l^2(\bz)\ni\xi\stackrel{W}{\rightarrow}F_\xi\in\mathcal H(\chi)
\end{equation}
It is clear that $W$ is isometric; indeed
$$\|F_\xi\|_{\chi}^2=\sum_{j\in\bz}|F_\xi(j,0)|^2=\sum_{j\in\bz}|\xi_j|^2,\quad(\xi\in l^2(\bz).$$
But it is also clear from Definition \ref{def2_5} that $W$ maps onto $\mathcal H(\chi)$; and so $W$ is a unitary isomorphism of the two Hilbert spaces.

The conclusion in our theorem now takes the form
\begin{equation}\label{eq30}
WU(g)W^*=\Ind_B^G(\chi)(g),\quad(g\in G),
\end{equation}
or equivalently
\begin{equation}\label{eq31}
WU(g)=\Ind_B^G(\chi)(g)W.
\end{equation}
The following computation proves the formula \eqref{eq31}. Let $\xi\in l^2(\bz)$, $g=(j,b)\in G, (k,c)\in G$; i.e., $k\in\bz$, $c\in B$. Then
$$(WU(g)\xi)(k,c)=\chi(c)(U(g)\xi)_k=\chi(c)(D_\chi(b)T_j\xi)_k=\chi(c)\chi(\alpha^k(b))\xi_{k+j}=\chi(\alpha^k(b)+c)\xi_{k+j}=$$
$$(W\xi)(k+j,\alpha^k(b)+c)=(W\xi)((k,c)g)=\Ind_B^G(\chi)(g)(W\xi)(k,c)$$
where we used \eqref{eq17} in Definition \ref{eq25} in the last step. 

Hence $WU(g)\xi=\Ind_B^G(\chi)(g)W\xi$ which is the desired conclusion \eqref{eq31} of the theorem.
\end{proof}

The next result shows that the family of unitary irreducible representations comes with a spectral dichotomy: If a given unitary irreducible representation, i.e., $U\in\Rep_{irr}(G,\mathcal H)$ has one discrete non-zero spectral component, then it is unitarily equivalent to $\Ind_B^G(\chi)$ for some $\chi\in K(=\hat B).$

We need a few preliminaries: By Stone's theorem, every $V\in\Rep(B,\mathcal H)$ decomposes as direct integral of one-dimensional representations, i.e., points $\chi$ in $K$. The decomposition has a discrete part and a continuous part; but either one may be zero. The discrete part $\mathcal H_p$ has the form 
\begin{equation}\label{eq4dz1}
\mathcal H_p(V)=\sum_{\chi}\mathcal H(\chi,V)
\end{equation}
(with each space $\mathcal H(\chi,V)$ being $B$-invariant)
where 
\begin{equation}\label{eq4dz2}
\mathcal H(\chi):=\mathcal H(\chi,V)=\{F\in\mathcal H\,|\, V(b)F=\chi(b)F, b\in B\}.
\end{equation}

\begin{corollary}\label{cor4dz}
Let $U\in\Rep_{irr}(G,\mathcal H)$, $U$ infinite dimensional, and consider the restriction $V(b):=U(0,b)$, $b\in B$. Suppose $\mathcal H_p(V)\neq\{0\}$. Then there is a $\chi\in K$ such that
\begin{equation}\label{eq4dz3}
U\cong \Ind_B^G(\chi),\mbox{ unitary equivalence};
\end{equation}
and 
\begin{equation}\label{eq4dz4}
\mathcal H={\sum_{j\in\bz}}^\oplus\mathcal H(\hat\alpha^j(\chi),V).
\end{equation}
\end{corollary}

\begin{proof}
Let $U$ be as stated in the corollary. If $\mathcal H_p\neq 0$, then there is a $\chi\in K$, and $F$ in $\mathcal H$ with $\|F\|=1$ such that
$$U(0,b)F=\chi(b)F,\quad(b\in B).$$
In other words, $\mathcal H(\chi,V)\neq0$. But then
\begin{equation}\label{eq4dz5}
U(0,b)U(j,0)F=\chi(\alpha^{-j}(b))U(j,0)F,\quad(j\in\bz,b\in B).
\end{equation}
This means that the right-hand side in \eqref{eq4dz4} is a closed non-zero $(U,G)$ invariant subpace in $\mathcal H$.
Since $U$ is irreducible, we must have equality in \eqref{eq4dz4}. Also, the space $\mathcal{H}(\chi,V)$ is one-dimensional, ortherwise the sum in the right-hand side of \eqref{eq4dz4} decomposes into two $(U,G)$ invariant subspaces in $\mathcal H$, contradicting again the irreducibility of $U$. 
We have $\ip{F}{U(0,b)U(j,0)F}=\ip{U(0,-b)F}{U(j,0)F}$ for all $j\in\bz$ and $b\in B$. This implies that
$$(\chi(b)-\chi(\alpha^j(b)))\ip{F}{U(j,0)F}=0,\quad(j\in\bz,b\in B).$$
Since the representation is infinite dimensional, $\chi$ is not periodic and $\ip{F}{U(j,0)F}=0$ for $j\neq0$, thus the sum in the right-hand side of \eqref{eq4dz4} is an orthogonal one.

Finally, composing \eqref{eq4dz5} with \eqref{eq25} and \eqref{eq29} from the proof of Theorem \ref{th3_2}, we conclude that the unitary equivalence assertion \eqref{eq4dz3} holds (just map $U(j,0)F$ into the canonical vectors $\delta_{-j}\in l^2(\bz)$ to construct the intertwining isomorphism), and the proof is completed.
\end{proof}
\begin{remark} 
We will classify the irreducible finite dimensional representations in Corollary \ref{cor4_3.5}
\end{remark}

\subsection{Periodic points and $\Ind_B^G(\chi)$}\label{periodic}
In this subsection we examine the commutant of the representations $\Ind_B^G(\chi)$. 
\begin{definition}\label{def3_3}
We say that a point $\chi\in K=\hat B$ is periodic of period $p$ if there is a $p\in\bn$ such that $\hat\alpha^p(\chi)=\chi$. We say that $p$ is {\it the period} if $\hat\alpha^k(\chi)\neq\chi$ for $1\leq k<p$, i.e., if $p$ is the first occurence of return to $\chi$.
\end{definition}

\begin{remark}\label{rem3_4}
For the cases $B=\bz_d[A^{-1}]$ in Example \ref{ex2_3} and \ref{ex2_4}, the set of periodic points in $K_A$ is countable. We give the details for Example \ref{ex2_3} but they extend {\it mutatis mutandis} to Example \ref{ex2_4}.

For fixed $p\in\bn$, the points $\chi=(z_0,z_1,z_2,\dots)$ in $K_2=\widehat{\bz[\frac12]}$ of period $p$ have the following form:

Let $z_0$ be a $(2^p-1)$'th root of $1$, i.e., a solution to $z^{2^p-1}=1$. Consider the finite word $\eta=(z_0,z_0^{2^{p-1}},\dots,z_0^2,z_0)$, and set $\chi=(\eta,\eta,\dots)$ infinite repetition of the finite word $\eta=\eta(z_0)$. Then $\chi\in K_2$ has period $p$; and conversely every $\chi$ in $K_2$ of period $p$ has this form.
\end{remark}

\begin{definition}\label{def3_5}
Returning to the general case, if some $\chi\in K=\hat B$ does not have a finite period, we say it is {\it aperiodic}. This happens iff the points $\{\hat\alpha^j(\chi)\,|\,j\in\bz\}$ are distinct, i.e., iff the mapping $\bz\ni j\mapsto\hat\alpha^j(\chi)\in\mathcal{O}(\chi)=:${\it the orbit of }$\chi$, is one-to-one.
\end{definition}

\begin{theorem}\label{th3_6}
(i) The representation $\Ind_B^G(\chi)$ is irreducible if and only if $\chi$ is aperiodic.

(ii) If some $\chi\in K$ has minimal period $p\in\bn$ (finite), then the commutant of $\Ind_B^G(\chi)$ is isomorphic to the abelian algebra of multiplication operators $\{f(z^p)\,|\,f\in L^\infty(\bt)\}$ where $\bt=\br/\bz$.
\end{theorem}

\begin{proof}
As noted in \eqref{eq29} from the proof of Theorem \ref{th3_2}, the Hilbert space which carries $\Ind_B^G(\chi)$ is $l^2(\bz)$. But it is convenient to make use of the isomorphism $l^2(\bz)\cong L^2(\bt)$ via Fourier series: If $\xi=(\xi_k)_{k\in\bz}\in l^2(\bz)$, we define 
\begin{equation}\label{eq3.16}
f_\xi(z)=\sum_{k\in\bz}\xi_k z^k,\quad(z\in\bt)
\end{equation}
If an operator $A:l^2(\bz)\rightarrow l^2(\bz)$ is in the commutant of some $\Ind_B^G(\chi)$, it follows from Lemma \ref{lem3_1} that it has the form $M_F$ for some $F\in L^\infty(\bt)$, i.e.,
\begin{equation}\label{eq3.17}
(Af_\xi)(z)=F(z)f_\xi(z),\quad(z\in\bt)
\end{equation}
If
\begin{equation}\label{eq3.18}
F(z)=\sum_{k\in\bz}\eta_kz^k
\end{equation}
then 
\begin{equation}\label{eq3.19}
(\chi(\alpha^n(b))-\chi(\alpha^k(b))\eta_{n-k}=0,\quad(n,k\in\bz).
\end{equation}
Conclusions: (i) If $\chi$ is aperiodic, then $\eta_k=0$ for all $k\in\bz\setminus\{0\}$, and we conclude that $A=\eta_0 I$ with $I$ denoting the identity operator in $L^2(\bz)$, or equivalently in $L^2(\bt)$. 

(ii) If $\chi$ has minimal period $p$, it follows that $\eta_k=0$ whenever $k\not\equiv0\mod p$ or, the possible non-zero terms have the form $\eta_{np}$, $n\in\bz$. Using \eqref{eq3.18}, we conclude that $F(z)=f(z^p)$ if we set $f(z):=\sum_{n\in\bz}\eta_{np}z^n$. This proves the conclusion in the theorem in both cases.
\end{proof}

\begin{definition}\label{def3_7}
Let $\chi\in K(=\hat B)$; then we say that the set 
\begin{equation}\label{eq3.20}
\mathcal O(\chi):=\{\hat\alpha^k(\chi)\,|\,k\in\bz\}
\end{equation} 
is the orbit of $\chi$.
\end{definition}

\begin{definition}\label{def3_8}
Consider two unitary representations $U_1$, and $U_2$, i.e., $U_i\in\Rep(G,\mathcal H_i)$, $i=1,2$. Set 
$$\mathcal L_G(U_1,U_2):=\{A:\mathcal H_1\rightarrow\mathcal H_2\,|\, A\mbox{ is bounded and linear, and }AU_1(g)=U_2(g)A,\, g\in G\}.$$
If $\mathcal L_G(U_1,U_2)=0$ we say that the two representations are {\it disjoint}. 
\end{definition}

\begin{corollary}\label{cor3_9}
Let $\chi_1,\chi_2\in K=\hat B$ and let $U_i:=\Ind_B^G(\chi_i)$, $i=1,2$ be the corresponding induced representations. Then $\mathcal L_G(U_1,U_2)\neq\{0\}$ if and only if $\mathcal O(\chi_1)\cap\mathcal O(\chi_2)\neq \emptyset$.
\end{corollary}

\begin{proof}
As we noted in the proof of Theorem \ref{th3_6}, an operator $A$ in $\mathcal L_G(U_1,U_2)$ must have the form $A=M_F$ for some $F\in L^\infty(\bt)$. For the Fourier expansion \eqref{eq3.18}, we get 
\begin{equation}\label{eq3.22}
(\chi_1(\alpha^n(b))-\chi_2(\alpha^k(b)))\eta_{n-k}=0,\quad(n,k\in\bz,b\in B).
\end{equation}
If $\mathcal O(\chi_1)\cap\mathcal O(\chi_2)=\emptyset$, for all $n,k\in\bz$, there exists $b\in B$ such that $\chi_1(\alpha^n(b))=\chi_2(\alpha^k(b))$, and conversely. The result now follows from \eqref{eq3.22}. 

\end{proof}

   The reader will notice that the ideas underlying our present discussion of Corolarry \ref{cor3_9} and Theorem \ref{th3_11}, below are very close to G. W. Mackey's view on unitary representations; see e.g., \cite{Mac76}.

\begin{definition}\label{def3_10}
We say that two representations $U_1$ and $U_2$ are {\it unitarily equivalent} iff there is a unitary isomorphism $W\in \mathcal L_G(U_1,U_2)$; i.e., $W:\mathcal H_1\rightarrow\mathcal H_2$, unitary (including ``onto'') such that
\begin{equation}\label{eq3.23}
WU_1(g)=U_2(g)W,\quad(g\in G).
\end{equation}
We shall also use the notation $U_1\cong U_2$; and we set 
$$\operatorname*{Class}(U_1)=\{U_2\,|\,U_2\cong U_1\}.$$
\end{definition}

\begin{theorem}\label{th3_11}
There is a natural bijection between the set of all orbits in $K=\hat B$ on one side, and on the other the set of all equivalence classes of induced representations of $G=B\rtimes_\alpha\bz$; with the bijection being implemented by
\begin{equation}\label{eq3.24}
K\ni\chi\mapsto\Ind_B^G(\chi)\in\Rep(G,l^2(\bz)).
\end{equation}
\end{theorem}

\begin{proof}
The details of the proof are essentially contained in the previous discussion. An essential point in the argument is that different points $\chi_1$ and $\chi_2$ in the same orbit are mapped into unitarily equivalent representations under \eqref{eq3.24}. To see this note that if $\chi_2=\hat\alpha^k\chi_1$ for some $k\in\bz$, then
\begin{equation}\label{eq3.25}
\Ind_B^G(\chi_2)=T_k\Ind_B^G(\chi_1)T_k^*.
\end{equation}
The formula \eqref{eq3.25} in turn follows from Lemma \ref{lem3_1}. This means that $T_k\in\mathcal L_G(\Ind_B^G(\chi_1),\Ind_B^G(\chi_2))$. 

Note that if the representations are not irreducible, there may be other intertwining operators.
\end{proof}

\begin{remark}\label{rem3_12}
Since the equivalence classes (referring to unitary
equivalence) of the induced representations are indexed by the set of orbits
in $K =\hat B$, one might naturally ask for a concrete and measurable cross
section for the $\hat\alpha$-orbits. Measurable cross sections play an essential role
in Mackey's theory of direct integrals, see e.g., \cite{Mac76}. However the
initial impetus for this theory came from the study of type $I$ groups. We
will see in section \ref{plancherel} that $G = B \rtimes_\alpha \bz$ is non-type $I$.
We show below in the case of $K_2 = \widehat{\bz[\frac12]}$ that in general we have
non-existence of measurable cross sections for the $\hat\alpha$-orbits.

 The fact that discrete semi-direct products of the form $B\rtimes_\alpha Z$ are not type $I$ was known to Mackey in the late 40's, early 50's (and in fact, this sparked his interest in ergodic theory).  In the spirit of \cite{Tho64}, we will prove this fact about the non-type $I$ status of $B\rtimes_\alpha\bz$ directly in Section 6. 
\end{remark}

\begin{proposition}\label{prop3_13} 
Let $\mu$ be the normalized Haar measure on the compact group $K=\hat B$. 
If the set of periodic points of $\hat\alpha$ has $\mu$-measure zero, then there is no measurable subset $M$ of $K$ such that $M$ intersects $\mathcal O(\chi)$ in exactly one point, for $\mu$-almost every $\chi\in K$.
\end{proposition}

\begin{proof}
Suppose {\it ad absurdum} that there is such a measurable subset $M$. Then, since we can eliminate the set of periodic points (having measure zero), we get that $\hat\alpha^k(M)\cap\hat\alpha^l(M)=\emptyset$ for all $k\neq l$ in $\bz$, up to measure zero, and $\bigcup_{k\in\bz} \hat\alpha^k(M)=K$. 
But, then 
$$1=\mu(K)=\mu(\bigcup_{k\in\bz}\hat\alpha^k(M))=\sum_{k\in\bz}\mu(\hat\alpha^k(M)).$$
But, from Lemma \ref{lem4_1} below, we see that $\mu(\hat\alpha^k(M))=\mu(M)$ for all $k\in\bz$. Either of the two possibilities for $\mu(M)$ yields a contradiction, either way $\mu(M) = 0$, or positive, contradicts the sum formula.
\end{proof}

\begin{remark}\label{rem3_14}
We define the equivalence relation $\sim$ on $K$ by $\chi_1\sim\chi_2$ iff $\mathcal O(\chi_1)=\mathcal O(\chi_2)$, and let $q:K\rightarrow K/_\sim$ be the quotient map. A section in $\sim$ is a map $m:K/_\sim\,\rightarrow K$ such that
$q\circ m=\mbox{id}_{K/_\sim}$. Proposition \ref{prop3_13} shows that, when the periodic points have $\mu$-measure zero, there are no measurable sections in $\sim$.
\end{remark}

\section{Haar measure}\label{haar}

For general references on elementary facts about ergodic automorphisms of compact abelian groups, invariance of Haar measure, etc, use for example \cite{Wal82, Pet83}.

In Definition \ref{def2_2}, we considered a dual automorphism $\hat\alpha\in\Aut(\hat B)$ arising from a fixed $\alpha\in\Aut(B)$ where $B$ is a given discrete abelian group. By $\hat B=:K$ we mean the compact dual group. As a compact abelian group, $K$ has a unique normalized Haar measure $\mu$, i.e., $\mu$ is a positive Borel measure on $K$ such that $\mu(K)=1$ and 
\begin{equation}\label{eq32}
\mu(E\chi)=\mu(E)\mbox{ for all } E\in\mathcal B(K)(=\mbox{ the Borel subsets in }K.)
\end{equation}
Here $E\chi:=\{\eta\chi\,|\,\eta\in E\}$, and $(\eta\chi)(b):=\eta(b)\chi(b)$, for all $b\in B$.

In general, if $\tau:K\rightarrow K$ is a measurable endomorphism, we set
\begin{equation}\label{eq33}
(\mu\circ\tau^{-1})(E)=\mu(\tau^{-1}(E)),\mbox{ for }E\in\mathcal B(K);
\end{equation}
and we note that $\mu\circ\tau^{-1}$ is again a Borel measure.

In this section we shall use that $\mu$ is automatically $\hat\alpha$-invariant, i.e., that 
\begin{equation}\label{eq34}
\mu\circ\hat\alpha^{-1}=\mu.
\end{equation}
In particular, this holds when $\hat\alpha=\hat\alpha_A$ is the automorphism induced by the matrix-examples, Example \ref{ex2_4}.

\begin{lemma}\label{lem4_1}
Let $\hat\alpha\in\Aut(K)$ be an automorphism of a compact abelian group $K$. Then the (normalized) Haar measure $\mu$ on $K$ is $\hat\alpha$-invariant.
\end{lemma}

\begin{proof}
It follows from the definitions that the measure $\mu\circ\hat\alpha^{-1}$ is translation invariant and normalized, i.e., that \eqref{eq32} holds. The conclusion now follows from uniqueness of the Haar measure. 
\end{proof}

\begin{corollary}\label{cor4_2}
Let $d\in\bn$, and let $A$ be a $d\times d$ matrix over $\bz$ such that $\det A\neq 0$. Let $\hat\alpha_A\in\Aut(K_A)$, $K_A:=\widehat{\bz_d[A^{-1}]}$, be the induced automorphism of Example \ref{ex2_4}, and let $\mu$ be the Haar measure on $K_A$. Then
\begin{equation}\label{eq35}
\mu\circ\hat\alpha_A^{-1}=\mu.
\end{equation}
\end{corollary}

\begin{definition}\label{def4_3}
Let $(K,\mu)$ be as above, and let $\hat\alpha\in\Aut(K)$ be given. We say that $\hat\alpha$ is {\it ergodic} if the only functions in $L^\infty(K)$ satisfying
\begin{equation}\label{eq37}
f\circ\hat\alpha=f
\end{equation}
are the constants, a.e., with respect to $\mu$. 
\end{definition}

\begin{lemma}\label{lem4_4}
Let $\hat\alpha_A\in\Aut(K_A)$ be the automorphism from Corollary \ref{cor4_2}. Assume that all the eigenvalues $\lambda$ of $A$ satisfy $|\lambda|>1$, i.e., that $A$ is an {\it expansive} matrix. Then $\hat\alpha_A$ is ergodic.
\end{lemma}

\begin{proof}
We will use a result from \cite{BrJo91} to the effect that $\hat\alpha_A$ is ergodic if and only if for all $n\in\bn$ we have the following implication: $b\in B$, $\alpha_A^n(b)=b$ implies $b=0$. 

But in view of \eqref{eq11}-\eqref{eq12}, the assertion $\alpha_A^n(b)=b$ takes the form 
\begin{equation}\label{eq38}
A^nb=b
\end{equation}
where the left-hand side in \eqref{eq38} is given by matrix-multiplication. Since $B=\bz_d[A^{-1}]\subset\br^d$, solutions in $B$ to \eqref{eq38} are column vectors in $\br^d$. But the number $1$ is not in the spectrum of $A^n$ for any $n\in\bn$; so $b=0$. This proves the ergodicity.

\end{proof}

\begin{theorem}\label{th4_5}
Let $\alpha\in\Aut(B)$ be as in Definition \ref{def2_1}, and let $\hat\alpha\in\Aut(\hat B)$ be the dual automorphism. Set $K=\hat B$ (the Pontryagin dual), and let $\mu$ be the Haar measure on $K$. Let $G=B\rtimes_\alpha\bz$.
\begin{enumerate}
\item For $g=(j,b)\in G$, and $\chi\in K$, and $f\in L^2(K,\mu)$ setting
\begin{equation}\label{eq39}
(U(g)f)(\chi)=\chi(b)f(\hat\alpha^j(\chi))
\end{equation}
we note that $U\in\Rep(G,L^2(K,\mu))$.
\item The representation $U$ in (i) is irreducible if and only if $\hat\alpha$ is ergodic. In particular, irreducibility holds for the $\hat\alpha_A$-representation of Lemma \ref{lem4_4}.
\end{enumerate}
\end{theorem}

\begin{proof}
(i) A direct computation using \eqref{eq39} shows that $U(g_1)U(g_2)=U(g_1g_2)$, $g_1,g_2\in G$, where the multiplication $g_1g_2$ refers to the product in $G$; see \eqref{eq1}. Since by Lemma \ref{lem4_1}, the Haar measure $\mu$ is preserved by $\hat\alpha$, it follows that $U$ is indeed a unitary representation.

(ii) Let $A$ be a bounded operator in $L^2(K,\mu)$ such that 
\begin{equation}\label{eq40}
AU(g)=U(g)A,\quad(g\in G).
\end{equation}
We say that $A$ is in the {\it commutant} of $U$. From \eqref{eq40}, we claim that $A$ must be a multiplication operator, i.e., that for some $f\in L^\infty(K)$, $A$ has the form
\begin{equation}\label{eq41}
(Au)(\chi)=f(\chi)u(\chi),\quad(u\in L^2(K,\mu),\chi\in K).
\end{equation}
Returning to \eqref{eq39} in the special case of $U(0,b)(\chi)=e_b(\chi):=\chi(b)$ we note that 
$U(0,c)e_b=e_{c+b}=e_ce_b$. This means that if $A$ satisfies \eqref{eq40} then $A$ must commute with all the multiplication operators for $f=e_b$, $b\in B$. But by Stone-Weierstrass, the linear combinations of $\{e_b\,|\,b\in B\}$ are dense in $C(K)$. Hence $A$ must itself be a multiplication operator, i.e., have the form \eqref{eq41}.

Since $A$ also commutes with $\{U(j,0)\,|\,j\in\bz\}$ we conclude that 
\begin{equation}\label{eq42}
f\circ\hat\alpha=f;
\end{equation}
in other words the commutant is the abelian algebra of all multiplication operators $M_f$ defined from $f\in L^\infty(K)$ of \eqref{eq42}. The results of (ii) now follows, see Definition \ref{def4_3}.
\end{proof}

Our next result yields the spectral type of the projection valued measure $P(\cdot)$ from Stone's formula applied to the restricted representation. 

We consider the representation $U\in\Rep(G,L^2(K,\mu))$ from \eqref{eq39} in Theorem \ref{th4_5}. Using the restriction $U|_B$ to the abelian subgroup $B$, we show that $U$ is disjoint from $\Ind_B^G(\chi)$ for all $\chi\in\hat B$.
To determine the projection valued measure $P(\cdot)$ on $K=\hat B$ in \eqref{eq14} for $V(b):=U(0,b)$, $b\in B$, we find
$$\ip{F}{P(dx)F}=\|P(dx)F\|^2,\mbox{ for all }F\in L^2(K,\mu).$$
We set $\nu_F(\cdot):=\|P(\cdot)F\|^2$. From \eqref{eq14}, we see that 
\begin{equation}\label{eq43}
\ip{F}{V(b)F}=\hat\nu_F(b),\quad(b\in B)
\end{equation}
where $\hat\nu_F$ denotes the Fourier transform of the measure $\nu_F$. 
\begin{lemma}\label{lem4_6}
Consider the representation $U$ in \eqref{eq39}, and its restriction to $B$, and let $P$ be the corresponding projection valued measure. Then 
\begin{equation}\label{eq44}
\nu_F(\cdot)=\|P(\cdot)F\|^2
\end{equation}
is absolutely continuous with respect to the Haar measure $\mu$ on $K$ for all $F\in L^2(K,\mu)$, and we have the following formula 
\begin{equation}\label{eq45}
\frac{d\nu_F}{d\mu}=|F|^2
\end{equation}
for the Radon-Nikodym derivative; i.e., we have $\nu_F\ll\mu$, and the Radon-Nikodym derivative is the $L^1$-function $|F|^2$ on $K$. 
\end{lemma}

\begin{proof}
We establish the identity \eqref{eq45} by checking that the Fourier transform applied to the two sides yields the same result. Specifically, we claim that 
\begin{equation}\label{eq46}
d\nu_F=|F|^2\,d\mu.
\end{equation}
For the Fourier transform, we have $(b\in B=\hat K)$:
$$\widehat{d\nu_F}(b)=\int_K\chi(b)\,d\nu_F(\chi)=\int_K\chi(b)\|P(d\chi)F\|_{L^2(\mu)}^2=\ip{F}{U(0,b)F}_{L^2(\mu)}=\int_K\cj{F(\chi)}\chi(b)F(\chi)\,d\mu(\chi)=$$
$$\int_K\chi(b)|F(\chi)|^2\,d\mu(\chi)=\widehat{(|F|^2\,d\mu)}(b).$$
Since $b\in B=\hat K$ is arbitrary, the desired formula \eqref{eq46} follows.
\end{proof}

The conclusions in Lemma \ref{lem4_6} states that the representation from Theorem \ref{th4_5} has continuous spectrum when restricted to $B$. But we must then know that the Haar measure $\mu$ on $K=\hat B$ does not have atoms. 

\begin{proposition}\label{prop4_7}
Let $K_A=\widehat{\bz_d[A^{-1}]}$ be the Pontryagin dual of $B_A=\bz_d[A^{-1}]$ from Example \ref{ex2_4}, and let $\mu=\mu_A$ denote the normalized Haar measure of $K_A$. Then $\mu_A$ is non-atomic, i.e., for all $\chi_0\in K_A$, 
\begin{equation}\label{eq47}
\mu_A(\{\chi_0\})=0,
\end{equation}
where $\{\chi_0\}$ denotes the singleton. 
\end{proposition}

\begin{proof}
Suppose $\chi_0$ is an atom. Then, since $\mu_A$ is invariant under
translations, every $\chi\in K_A$ is an atom and
$\mu_A(\{\chi\})=\mu_A(\chi\chi_0^{-1}\{\chi_0\})=\mu_A(\{\chi_0\})$.

Since $K_A$ is compact and $\mu_A(K_A)=1$, it follows that for any finite
subset of $K_A$, $$1=\mu_A(K_A)\geq \sum_{\chi\in F}\mu_A(\{\chi\})= |F| \mu_A(\{\chi_0\})$$
and as $K_A$ is infinite, this implies that $\mu_A(\{\chi_0\})=0$, a contradiction.
\end{proof}

This is a statement of a more general result about compact connected abelian groups. Note the dual of $K_A$
has no torsion and is countable discrete abelian so $K_A$ is a connected compact abelian group; see \cite{HeRo63}.

\begin{corollary}\label{cor4_8}
Let $A$, and $(K_A,\mu_A)$ be as in the statement of Proposition \ref{prop4_7}. Let $\chi_0\in K_A$. Then 
$$\{F\in L^2(K_A,\mu_A)\,|\,U(0,b)F=\chi_0(b)F,\,b\in B_A\}=\{0\};$$
i.e., the representation 
\begin{equation}\label{eq50}
(U(0,b)F)(\chi)=\chi(b)F(\chi),\quad(\chi\in K_A)
\end{equation}
in $L^2(K_A,\mu_A)$ has no point-spectrum. 
\end{corollary}

This is a very general fact about ergodic automorphisms on compact connected abelian groups. They never have point spectrum; see \cite{Wal82,Pet83}.
\begin{proof}
Suppose for some $\chi_0\in K_A$, and $F_0\in L^2(K_A,\mu_A)$ we had 
$$U(0,b)F_0=\chi_0(b)F_0,\quad(b\in B_A).$$
Then \eqref{eq50} shows that $F_0$ is a constant times the indicator function of $\chi_0$. But the indicator functions of points in $K_A$ have $L^2(\mu_A)$-norm equal to $0$.  
\end{proof}

\section{The Plancherel formula for $G$}\label{plancherel}

       In sections \ref{spectral} and \ref{ind} we examined the monomial representations of the discrete semidirect product group $G = B \rtimes_\alpha \bz$. We recall that the starting point is a given discrete abelian group $B$, and a fixed automorphism $\alpha$ of $B$. By a {\it monomial representation} of $G$ we mean a representation of $G$ which is induced from a one-dimensional representation of $B$. Since the one-dimensional representation of $B$ are the points $\chi$ in the compact Pontryagin dual group $K:= \hat B$, the monomial representations of $G$ have the form $\Ind_B^G(\chi)$, and they are indexed by points $\chi$ in $K$.

       Note that since the group $G = B \rtimes_\alpha\bz$ is discrete, its (right) regular representation $R$ simply acts by right translations on the sequence $l^2$ Hilbert space $l^2(G)$. In this section we prove that $R$ is the direct integral of the induced representations $\Ind_B^G(\chi)$ with the integration being with respect to the Haar measure $\mu$ on $K$. This means that the particular continuous representation in section \ref{haar} does not contribute to the Plancherel formula for $G$. Or stated differently, only the representations of $G$ whose restriction to $B$ have point-spectrum contribute to the Plancherel formula for $G$.

       It is interesting to compare this result to the theorem in \cite{LPT01}. The authors of \cite{LPT01}, see section \ref{spectral} above, show that for affixed dilation matrix $A$, the direct integral of certain ``thin'' subsets, the wavelet sets, in $K$ ``add up'' to the $A$-wavelet representation in $L^2(\br^d)$; see also Example \ref{ex2_4}.

Specifically, if $d\in\bn$ and $A$ is a $d\times d$ matrix as in Example \ref{ex2_4}, we consider the Hilbert space $L^2(\br^d)$ relative to the usual $d$-dimensional Lebesgue measure on $\br^d$. 

We define the wavelet representation of $G_A=B_A\rtimes_{\alpha_A}\bz$ on $L^2(\br^d)$ as follows  (see also \cite{MaVa00}):
\begin{equation}\label{eq51}
(U_w(j,A^{-n}k)f)(x)=|\det A|^{j/2}f\left(A^j(x-A^{-n}k)\right),\quad(f\in L^2(\br^d),x\in\br^d,k\in\bz^d,j\in\bz,n\geq 0).
\end{equation}
As in Example \ref{ex2_4}
\begin{equation}\label{eq52}
B_A:=\bz_d[A^{-1}]=\bigcup_{n\geq 0}A^{-n}\bz^d\subset\br^d.
\end{equation}
Points $x\in\br^d$, and $k\in\bz^d$ are viewed as column vectors. 
\par
Note that the wavelet representation is obtained from translation and dilation operators: we have 
\begin{equation}\label{eqd52_1}
(T_{A^{-n}k}f)(x)=(U_w(0,A^{-n}k)f)(x)=f(x-A^{-n}k),\quad (Df)(x)=(U_w(1,0)f)(x)=|\det A|^{1/2}f(Ax),
\end{equation}
and 
\begin{equation}\label{eqd52_2}
U_w(j,A^{-n}k)=T_{A^{-n}k}D^j,\quad(j\in\bz,A^{-n}k\in\bz_D[A^{-1}]).
\end{equation}

We say that $U_w$ is the $A$-wavelet representation.

\begin{definition}\label{def5_1}
A {\it wavelet set} is a measurable subset $E\subset\br^d$ whcih satisfies the following four conditions:
\begin{enumerate}
\item $\bigcup_{k\in\bz}(E+k)=\br^d$;
\item $E\cap(E+k)=\empty$ if $k\in\bz^d\setminus\{0\}$;
\item $\bigcup_{j\in\bz}A^jE=\br^d$;
\item $E\cap A^jE=\empty$ if $j\in\bz\setminus\{0\}$.
\end{enumerate}
\end{definition}

\begin{remark}\label{rem5_2}
Note that the four conditions (i)--(iv) are occasionally stated to hold only up to $\br^d$-Lebesgue measure zero. But since the operations in (i) and (iii) are countably discrete, a given set $E$ which satisfies the conditions modulo $d$-Lebesgue measure zero, may be modified so that (i)--(iv) hold everywhere, the modifications involving only changes on subsets of measure zero.
\end{remark}

The theorem in \cite{LPT01} states that a subset $E\subset\br^d$ is a $A$-wavelet set if and only if it supports and orthogonal direct integral decomposition for $U_w$ of \eqref{eq51}; i.e., iff
\begin{equation}\label{eq53}
U_w={\int_E}^\oplus\Ind_{B_A}^{G_A}(\chi_t)\,dt
\end{equation}
where ``$dt$'' is the $d$-dimensional Lebesgue measure supported on the set $E$. Actually, the theorem in \cite{LPT01} shows that one can perform the decomposition in \eqref{eq53} even if the set $E$ tiles $\br^d$ only by dilations by $A^T$. The tranlation tiling is not needed for the decomposition and was only used elsewhere. Moreover, for $t\in E(\subset \br^d)$ we use the notation $\chi_t\in K_A=\widehat{B_A}=\widehat{\bz_d[A^{-1}]}$. Specifically, we use the setting in Example \ref{ex2_4}. Recall that the range of the $\br^d$-embedding in $K_A$ is dense, i.e., $\br^d\ni t\mapsto\chi_t\in K_A$ coming from dualizing $B_A=\bigcup_{k\geq0}A^{-k}\bz^d\subset\br^d$.

\begin{theorem}\label{th5_3}
Let $B$ be a discrete abelian group, $\alpha\in\Aut(B)$, and $G=B\rtimes_\alpha\bz$. Let $R:G\rightarrow\mathcal B(l^2(G))$ (in fact unitary operators on $l^2(G)$) be the regular representation 
\begin{equation}\label{eq5.3}
(R(g_1)F)(g_2)=F(g_2g_1),\quad(g_1,g_2\in G, F\in l^2(G)).
\end{equation}
For $\chi\in K:=\hat B$ set 
\begin{equation}\label{eq5.4}
U_\chi:=\Ind_B^G(\chi);
\end{equation}
see section \ref{ind} above. Let $\mu$ be the normalized Haar measure of the compact group $K$. Then we have the following orthogonal direct integral representation 
\begin{equation}\label{eq5.5}
R={\int_K}^\oplus U_\chi\,d\mu(\chi).
\end{equation}
\end{theorem}

\begin{proof}
There are several parts to the formula \eqref{eq5.5}. First recall the regular representation $R$ acts on the Hilbert space $l^2(G)$ while each $U_\chi=\Ind_B^G(\chi)$ acts on $l^2(\bz)\cong L^2(\bt)$ as described in Theorem \ref{th3_2}. So part of the conclusion in \eqref{eq5.5} is the assertion 
\begin{equation}\label{eq5.6}
l^2(G)\cong{\int_K}^\oplus l^2(\bz)_\chi\,d\mu(\chi)
\end{equation}
where \eqref{eq5.6} is really an integral transform applied to elements $F\in l^2(G)$, i.e., $F:G\rightarrow\bc$ such that 
\begin{equation}\label{eq5.7}
\|F\|_{l^2(G)}^2:=\sum_{g\in G}|F(g)|^2.
\end{equation}

{\bf The transform.} For $\chi\in K$, and $F\in l^2(G)$, set 
\begin{equation}\label{eq5.8}
f_\chi(j,b)=\sum_{c\in B}\cj{\chi(c)}F(j,b+c),
\end{equation}
where in \eqref{eq5.8} we must for the moment restrict to $F$ of finite support. 

{\bf Two assertions:} For points $(j,b)\in G$
\begin{enumerate}
\item $f_\chi(j,b)=\chi(b)f_\chi(j,0)$ for all $j\in\bz, b\in B$;
\item 
$$\int_K\|f_\chi\|_{\mathcal H(\chi)}^2\,d\mu(\chi)=\|F\|_{l^2(G)}^2.$$
\end{enumerate}

First note that by \eqref{eq15} in Definition \ref{def2_5}, $f_\chi\in\mathcal H(\chi)=$ the space of $\Ind_B^G(\chi)$, for all $\chi\in K$.

{\it Proof of (i).}
$$f_\chi(j,b)=\sum_{c\in B}\cj{\chi(c)}F(j,b+c)=\sum_{c\in B}\cj{\chi(c-b)}F(j,c)=\chi(b)\sum_{c\in B}\cj{\chi(c)}F(j,c)=\chi(b)f_\chi(j,0)$$

{\it Proof of (ii).}
$$\int_K\|f_\chi\|_{\mathcal H(\chi)}^2\,d\mu(\chi)\stackrel{\mbox{\eqref{eq16}}}{=}\int_K\sum_{j\in\bz}|f_\chi(j,0)|^2\,d\mu(\chi)=\int_K\sum_{j\in\bz}\left|\sum_{b\in B}\cj{\chi(b)}F(j,b)\right|^2\,d\mu(\chi)=$$
$$\sum_{j\in\bz}\int_K\left|\sum_{b\in B}\cj{e_b(\chi)}F(j,b)\right|^2\,d\mu(\chi)\stackrel{\mbox{Plancherel}}{=}\sum_{j\in\bz}\sum_{b\in B}|F(j,b)|^2=\sum_{g\in G}|F(g)|^2=\|F\|_{l^2(G)}^2.$$
Thus the transform is isometric.

It remains to prove that the transform:
\begin{equation}\label{eq5.12}
l^2(G)\ni F\mapsto(f_\chi)\in{\int_K}^\oplus\mathcal{H}(\chi)\,d\mu(\chi)
\end{equation}
is onto. We will do this by exibiting an inverse.

First consider $f_\chi$ as in \eqref{eq5.8}. We claim that
\begin{equation}\label{eq5.13a}
\int_Kf_\chi(j,b)\,d\mu(\chi)=F(j,b),\quad((j,b)\in G).
\end{equation}
Indeed, in view of the isometric property of \eqref{eq5.12}, we may exchange integral and summation in \eqref{eq5.13a}. We get
$$\int_Kf_\chi(j,b)\,d\mu(\chi)=\sum_{c\in B}\int_K\cj{\chi(c)}F(j,b+c)\,d\mu(\chi)=F(j,b),$$
since $\int_K\cj{\chi(c)}\,d\mu(\chi)=\left\{\begin{array}{cc}1&\mbox{if }c=0\\
0&\mbox{if }c\neq0\end{array}\right.$.

The most general element in the direct-integral Hilbert space on the right-hand side in \eqref{eq5.12} is a measurable field $\varphi:K\times G\rightarrow\bc$ such that
\begin{equation}\label{eq5.13}
\varphi(\chi,bg)=\chi(b)\varphi(\chi,g),\quad(\chi\in K,b\in B,g\in G).\
\end{equation}
with 
\begin{equation}\label{eq5.14}
\int_K\|\varphi(\chi,\cdot)\|_{\mathcal{H}(\chi)}^2<\infty.
\end{equation}
If a meaurable field $\varphi$ is given subject to \eqref{eq5.13}--\eqref{eq5.14}, we may define
$$F(g)=\int_K\varphi(\chi,g)\,d\mu(\chi).$$
The previous computation shows that $F\in l^2(G)$, and that $F\mapsto f_\chi=\varphi(\chi,\cdot)$.

That the operator $F\mapsto (f_\chi)$ in \eqref{eq5.8} {\it intertwines} the respective representations amounts to the following identity

(iii) For $F\in l^2(G)$, and $g\in G$, we have
\begin{equation}\label{eq5.15}
(R(g)F)_\chi=\Ind_B^G(g)f_\chi,\quad(\chi\in K,g\in G).
\end{equation}
To prove \eqref{eq5.15}, let $g=(j,b)$, and evaluate the two sides in \eqref{eq5.15} at points $(k,c)\in G$; i.e., $b,c\in B$, and $j,k\in\bz$. Then
$$(R(g)F)_\chi(k,c)=\sum_{a\in B}\cj{\chi(a)}(R(g)F)(k,c+a)=\sum_{a\in B}\cj{\chi(a)}F(k+j,\alpha^k(b)+c+a)=\chi(\alpha^k(b)+c)f_\chi(k+j,0)=$$
$$f_\chi(k+j,\alpha^k(b)+c)=f_\chi((k,c)(j,b))=\Ind_B^G(g)f_\chi(k,c)$$
which is the desired formula.
\end{proof}

\begin{remark}
This is a folklore result and Mackey, Fell et al used this type of decomposition often (induction in stages result). The regular representation is just $\Ind_{\{0\}}^G(1)$ which is by induction in stages $\Ind_B^G(\Ind_{\{0\}}^B(1))$  and by Pontryagin duality $\Ind_{\{0\}}^B(1)=\int_{\hat B}^\oplus (\chi)d\,\chi$. Now use the fact that direct integrals and inducing commute, one of Mackey's favorite tricks.

\end{remark}

\begin{remark}\label{rem5_4}
In the last step of the computation for (ii) in the proof we replaced a $(K,\mu)$ integral with a $b\in B$ summation. This is based on the Pontryagin duality of $l^2(B)\cong L^2(K,\mu)$. In this duality, we have the following ONB in $L^2(K,\mu)$: $\{e_b\,|\,b\in B\}$ where $e_b:K\rightarrow\bc$ is defined by $e_b(\chi)=\chi(b)$ for all $b\in B, \chi\in K$. 

An important point of Pontryagin duality is that if $K=\hat B$ in the category of locally compact abelian groups, then $\hat K\cong B$ with a natural isomorphism.

The fact that $\{e_b\,|\,b\in B\}$ is an ONB follows from general Pontryagin duality, see e.g. \cite{Rud62}.
\end{remark}

\begin{remark}\label{rem5_5} The Baumslag-Solitar group, in the present context in the form $G = B \rtimes_\alpha \bz$ is an ICC group (meaning that its set of conjugacy classes is infinite, see e.g., \cite{Mac76}). By von Neumann's theory, this means that its right regular representation, i.e., $R$ in \eqref{eq5.3} in Theorem \ref{th5_3} will generate a von Neumann algebra factor of type $II_1$.

This does not contradict our direct integral decomposition \eqref{eq5.5} for $R$ into a direct integral of the family of monomial representations. Naturally the irreducible monomial representations give factors of type $I$. But the direct integral should be of type $II$.  

This is not a contradiction in view of Remark \ref{rem3_12}, i.e., non-existence of measurable cross sections in $K$.

Moreover, even though we have a direct integral decomposition \eqref{eq5.5}, this is not done over the center of the algebra. We have plenty of multiplicity in \eqref{eq5.5}, i.e., repetition of equivalent irreducible representations: All the representations from $\chi$ in the same orbit are equivalent by Theorem \ref{th3_11}. Specifically, every point in $K$ has an orbit under $\hat\alpha$ , and the representations corresponding to points in this orbit are equivalent.

Therefore each point will come with a set of operators that intertwine these irreducible representations along the orbit. Integrating will generate a big $II_1$ commutant.

        This is an important distinction between \eqref{eq53} and \eqref{eq5.5}. The first multiplicity free, and the second, far from it!

As noted, in the formula \eqref{eq5.5} there is a lot of multiplicity.

        The important point (Theorem \ref{th3_11}) is that the unitary equivalence classes of the representations $\Ind_B^G(\chi)$ are indexed by the $\hat\alpha$ orbits. There is a countable set of different $\chi$'s in the same orbit in \eqref{eq5.5}, so obviously the commutant corresponding to a fixed orbit $\mathcal O(\chi)$ is quite big. 
\end{remark}

\subsection{The ICC-condition}
The ICC-condition may be illustrated more clearly by use of Example \ref{ex2_3}. The issue is the set of conjugacy classes in the group $G=B\rtimes_\alpha\bz$. For $g_1,g_2\in G$ we say that $g_1\sim g_2$ (conjugacy) iff there is a $g_3\in G$ such that $g_2=g_3g_1g_3^{-1}$. The conjugacy class of $g\in G$ is denote $\tilde g$.

In Example \ref{ex2_3}, there is a natural system of bijections between the following three sets:
(i) Conjugacy classes $\{\tilde b\,|\, b\in B\}$; (ii) $\alpha$-orbits in $B$, i.e., $\operatorname*{Orb}(b)=\{\alpha^{j}b\,|\,j\in\bz\}$; and (iii) the set of odd integers. 
\begin{proof}
The assertion (i)$\Leftrightarrow$(ii) holds more generally; and follows from this: If $j\in\bz$, $b,c\in B$, then with $g=(j,c)$, we have 
\begin{equation}\label{eq6dz}
g^jbg^{-j}=\alpha^j(b).
\end{equation}
The remaining argument (ii)$\Leftrightarrow$(iii) uses the representation $B=\bz[\frac12]$ from Example \ref{ex2_3}.

For all $b\in B$ there is a unique $j\in\bz$ such that $2^jb\in\bz$, but $2^kb\not\in\bz$ if $k<j$. Set $\operatorname*{Ord}(b):=2^jb$. Then it is easy to see that the mapping $b\mapsto\operatorname*{Ord}(b)$ induces a bijection between the two sets in (ii) and (iii). 
\end{proof}

We mentioned that the regular representation $R$ in Theorem \ref{th5_3} generates a type $II_1$ factor von Neumann algebra of operators on $l^2(G)$. The trace $\tau(\cdot)$ on this factor is 
$$\tau(\cdot):=\ip{\delta_e}{\cdot\delta_e}_{l^2(G)}$$
where $e=(0,0)\in G$ is the neutral element.

\section{Finite dimensional representations}\label{finite}

We saw that the induced representation $\Ind_B^G(\chi)$ is reducible if $\chi\in K=\hat B$ has finite period. On the other hand it is still infinite-dimensional. The finite-dimensional representations are not induced from $B$ to $G=B\rtimes_\alpha\bz$. 

Consider $\chi$ of minimal period $p$, that is $\chi\in K$, and suppose $\hat\alpha^p\chi=\chi$, $\hat\alpha^k\chi\neq\chi$ if $1\leq k<p$. 
\begin{definition}
Set $B(\chi):=B/\{b\in B\,|\,\chi(b)=1\}$. Then $\alpha$ induces an action of $\bz_p=\bz/p\bz$ on $B(\chi)$. Now set
$G(\chi):=B(\chi)\rtimes_\alpha\bz_p$
\end{definition}

\begin{theorem}\label{th7_1}
The induced representation $\Ind_{B(\chi)}^{G(\chi)}(\chi)$ is $p$-dimensional and irreducible. 
\end{theorem}

\begin{proof}
The argument follows closely the one for induction $B\rightarrow G$ in section \ref{ind}, so we will only sketch the details. The important point is that the quotient qroup $B(\chi)\backslash G(\chi)$ is now a copy of $\bz_p=\bz/p\bz$. Hence the formula \eqref{eq25} for the $B\rightarrow G$ case modifies as follows
$$\Ind_{B(\chi)}^{G(\chi)}(\chi)=D_\chi^{(p)}(b)T_j^{(p)}$$
where
$$D_\chi^{(p)}(b)=\left(\begin{array}{cccc}
\chi(b)&0&\dots&0\\
0&\chi(\alpha(b))&\dots&0\\
\vdots& &\ddots&\vdots\\
0&\dots & &\chi(\alpha^{p-1}(b))\end{array}\right)$$
and $T_j^{(p)}:=(T^{(p)})^j$ with 
$$T^{(p)}:=\left(\begin{array}{cccccc}
0&1&0&\dots&0&0\\
0&0&1&\dots&0&0\\
\vdots&\vdots& & & &\vdots\\
0&0&0&\dots&1&0\\
0&0&0&\dots&0&1\\
1&0&0&\dots&0&0\end{array}\right)$$
The proof of irreducibility is moddeled on the argument for the proof of Theorem \ref{th3_6}(i) above.

\end{proof}

\begin{remark}\label{rem7_2}
Note that for $G=B\rtimes_\alpha\bz$ the formal Fr{\" o}benius reciprocity prediction breaks down and in fact:

\begin{theorem}\label{th7_3}
Let $\chi\in K$ be an element of finite period $p$, i.e., $\hat\alpha^p\chi=\chi$, $\hat\alpha^k\chi\neq\chi$, for $1\leq k<p$. Let $U_p^{(\chi)}$ be the finite-dimensional irreducible representation of $G$ given in Theorem \ref{th7_1}. Then $\mathcal L_G(U_p^{(\chi)},\Ind_B^G(\chi))=0$. 
\end{theorem}

\begin{proof}
We will write out the details only for $p=3$ to simplify notation. The general argument is the same. Recall
$$T^{(3)}=\left(\begin{array}{ccc}
0&1&0\\
0&0&1\\
1&0&0\end{array}\right),\quad D_\chi(b)=\left(\begin{array}{ccc}\chi(b)&0&0\\
0&\chi(\alpha(b))&0\\
0&0&\chi(\alpha^2(b))\end{array}\right)$$
\begin{equation}\label{eq10.3}
U_3^{(\chi)}(j,b)=D_\chi^{(3)}(b){T^{(3)}}^j
\end{equation}
while 
\begin{equation}\label{eq10.4}
(\Ind_B^G(\chi)_{(j,b)}\xi)_k=\chi(\alpha^k(b))\xi_{k+j},\quad((j,b)\in G, k\in\bz,\xi\in l^2(\bz)).
\end{equation}
Let $W\in \mathcal L_G(\U_p^{(\chi)},\Ind_B^G(\chi))$. Let $u_0,u_1,u_2$ be the canonical basis in $\mathcal H(U_p^{(\chi)})=\bc^3$.

Working mod $3$
\begin{equation}\label{eq10.6}
U_p^{(\chi)}(j,b)u_k=\chi(\alpha^{k+2j}(b))u_{k-j(\mod 3)},\quad((j,b)\in G,k\in\{0,1,2\}\cong\bz/3\bz).
\end{equation}
Set $Wu_k=\xi^{(k)}\in l^2(\bz)$, i.e., $\xi^{(k)}=(\xi_s^{(k)})_{s\in\bz}\in l^2(\bz)$, $\|\xi^{(k)}\|^2=\sum_{s\in\bz}|\xi_s^{(k)}|^2<\infty$.
Using \eqref{eq10.4}--\eqref{eq10.6}
$$WU_p^{(\chi)}(j,b)u_k=\Ind_B^G(\chi)(j,b)Wu_k,\quad((j,b)\in G,k\in\{0,1,2\})$$
so $$\chi(\alpha^{k-j}(b))\xi_s^{(k-j)_3}=\chi(\alpha^s(b))\xi_{s+j}^{(k)},\quad(s,j\in\bz).$$
Now set $j=3t\in 3\bz$, and we get $\chi(\alpha^k(b))\xi_s^{(k)}=\chi(\alpha_s(b))\xi_{s+3t}^{(k)}$ and 
$|\xi_s^{(k)}|=|\xi_{s+3t}^{(k)}|$ for all $s,t\in\bz$. 

Since $\xi^{(k)}\in l^2(\bz)$, we conclude that $\xi^{(k)}=0$ in $l^2(\bz)$; reason: $\lim_{t\rightarrow\infty}\xi_{s+3t}^{(k)}=0$.
\end{proof}

\end{remark}

\begin{corollary}\label{cor4_3.5}
If $U$ is an irreducible  finite dimensional representation of $G=B\rtimes_\alpha\bz$ then there exists a periodic element $\alpha\in K=\hat B$ of period $p$, and a $z_0\in\bt$ such that $U$ is unitarily equivalent to the representation
$$U_{\chi,z_0}(j,k)=D_\chi^{(p)}(b)T_{j,z_0}^{(p)},\quad(j\in\bz,b\in B)$$ where
$$D_\chi^{(p)}(b)=\left(\begin{array}{cccc}
\chi(b)&0&\dots&0\\
0&\chi(\alpha(b))&\dots&0\\
\vdots& &\ddots&\vdots\\
0&\dots & &\chi(\alpha^{p-1}(b))\end{array}\right), T_{j,z_0}^{(p)}:=(T_{z_0}^{(p)})^j,\mbox{ with }
T_{z_0}^{(p)}:=\left(\begin{array}{cccccc}
0&1&0&\dots&0&0\\
0&0&1&\dots&0&0\\
\vdots&\vdots& & & &\vdots\\
0&0&0&\dots&1&0\\
0&0&0&\dots&0&z_0\\
1&0&0&\dots&0&0\end{array}\right)$$
 The representations $U_{\chi,z_0}$ are disjoint for different pairs $(\chi,z_0)$.
\end{corollary}

\begin{proof}
The argument follows closely that of Corollary \ref{cor4dz}. Let $U\in\Rep_{irr}(G,\mathcal H)$, and suppose $\dim\mathcal H<\infty$. Then $\mathcal H_p\neq 0$ and by \eqref{eq4dz4} there are $\chi\in K$ and $F_0\in\mathcal H\setminus \{0\}$ such that $U(0,b)F_0=\chi(b)F_0$ for all $b\in B$. 

Assume $\|F_0\|=1$ and set $v_j:=U(-j,0)F_0$, $j\in\bz$. Then $U(0,b)v_j=\chi(\alpha^j(b))v_j$ for all $j\in\bz$, $b\in B$. From $\ip{v_j}{U(0,b)v_k}=\ip{U(0,-b)v_j}{v_k}$, an easy calculation yields
\begin{equation}\label{eqirr}
(\chi(\alpha^j(b))-\chi(\alpha^k(b)))\ip{v_j}{v_k}=0,\quad(j,k\in\bz).
\end{equation}
Since $\dim\mathcal H<\infty$, we conclude that $\chi$ has finite period. Let $p$ be the minimal period of $\chi$.
\par
From \eqref{eqirr}, it follows that $\ip{v_j}{v_k}=0$ if $j\not\equiv k\mod p$. 

Suppose we have the decomposition \begin{equation}\label{eq6star}
v_p=U(1,0)v_{p-1}=z_0v_0+w\end{equation} with $z_0\in\bc$, and $w\in\mathcal H\ominus\{v_0\}$.

We use that $\dim\mathcal H(\chi)=1$ (because the representation is irreducible, see also the proof of Corollary \ref{cor4dz}), so $\mathcal H(\chi)=\bc v_0$. Now apply $V_b$ $(b\in B)$ to \eqref{eq6star}: $\chi(b)v_p=V_bv_p=z_0\chi(b)v_0+V_bw$ for all $b\in B$. This implies $\chi(b)(v_p-z_0v_0)=V_bw$ so $w\in\mathcal H(\chi)=\bc v_0$. Hence $\mathcal H_p:=\mbox{span}\{v_0,v_1,\dots,v_{p-1}\}\subset\mathcal H$ is a $(U,G)$-invariant subspace, and $\{v_j\}_{j=0}^{p-1}$ is an ONB in $\mathcal H_p$. Since $(U,G)$ is irreducible it follows that $\mathcal H_p=\mathcal H$. Mapping the vectors $v_j$ into the canonical vectors $\delta_j$ in $\bc^p$ proves the first part of the corollary.

For disjointness, note that if $\chi_1\neq\chi_2$ then the $D_{\chi}^{(p)}$ parts of the representations cannot be equivalent. If $z_0\neq z_0'$ then, since the spectrum of the operator $T_{z_0}^{(p)}$ is the $p$-th roots of $z_0$, the $T_{z_0}^{(p)}$ parts of the two representations cannot be equivalent. Since the representations are irreducible, they must be disjoint.

\end{proof}

We should add that ecause $C^*(B\rtimes_A \bz)$ can be viewed as a
transformation group $C^*$-algebra $C(K_A)\rtimes\bz$, much facts about ideal structure and
representation theory can be gleaned from work of D. Williams in late 1970's.

\section{Dilations}

  One of the basic ideas underlying the multiresolution method in wavelet theory, and more generally in the analysis of fractals, see e.g., \cite{Dau95} and \cite{Jor06} is to select an initial Hilbert space, often as a subspace of an ambient Hilbert space. The ambient Hilbert space may be $L^2(\br^d)$, see e.g, Remark \ref{rem5_2}, and \cite{LPT01}. If the ambient Hilbert space is not given at the outset as is typically the case for fractals, e.g., \cite{DuJo06}, we must construct it from the initial data.

But in the initial space the scaling operator is merely an isometry, say $u$. The relation which governs the two wavelet operators in dyadic case is then $u t = t^2 u$, but after it is dilated (extended) via an isometric embedding, it may be realized with a version of u which will now be unitary, and we end up with a unitary representation of a group of the form $G = B \rtimes_\alpha \bz$ of the kind studied in sections \ref{spectral} to \ref{plancherel}. The isometric embedding must of course intertwine the two operator systems, and the dilated system must be minimal. In this section, we make this precise, and prove the properties which are asserted.

Let $H$ be a subgroup of $B$ with the property that
\begin{equation}\label{eqdil1}
\alpha(H)\subset H,\quad \bigcup_{n\geq0}\alpha^{-n}(H)=B.
\end{equation}
The inclusion $i:H\hookrightarrow B$ dualizes to a morphism $\theta_0:\hat B\rightarrow\hat H$.
\par
Let
\begin{equation}\label{eqdil2}
\mu_H:=\mu\circ\theta_0^{-1}.
\end{equation}
It is easy to see then that $\mu_H$ is a translation invariant measure, and $\mu_H(H)=1$, hence $\mu_H$ is the normalized Haar measure on $H$. 

Since $\alpha(H)\subset H$, the automorphism $\alpha$ on $B$ restricts to an injective endomorphism $\alpha_H$ of $H$, and it has dual morphism $\hat\alpha_H:\hat H\rightarrow\hat H$, which satsifies:
\begin{equation}\label{eqdil3}
\theta_0\circ\hat\alpha=\hat\alpha_H\circ\theta_0
\end{equation}

We define the operators $u$ and $t_h$, $h\in H$ on $L^2(\hat H,\mu_H)$ by
\begin{equation}\label{eqdil4}
(uf)(z)=f(\hat\alpha_H(z)),\quad (t_hf)(z)=z(h)f(z),\quad(f\in L^2(\hat H,\mu_H),h\in H,z\in\hat H).
\end{equation}

\begin{example}\label{example}
Consider the set-up of Example \ref{ex2_4}, i.e., $B=\bz_d[A^{-1}]$ and $\alpha(b)=Ab$ for $b\in B$, where $A$ is a $d\times d$ expansive integer matrix. 

Let $H=\bz^d$ as a subgroup of $\bz_d[A^{-1}]$. Then since $A$ is an integer matrix $\alpha(H)\subset H$ and by construction \eqref{eqdil1} is satisfied. 

The dual of the group $H=\bz^d$ is the $d$-torus $\hat H=\bt^d$, with the duality given by
$$\ip{(k_1,\dots,k_d)}{(e^{2\pi i x_1},\dots,e^{2\pi i x_d})}=e^{2\pi i\sum_{l=1}^dk_jx_j},\quad(k_j\in\bz,e^{2\pi ix_j}\in\bt).$$
The morphism $\hat\alpha_H$ on $\bt^d$ is then
$$\hat\alpha_H(e^{2\pi ix_1},\dots,e^{2\pi i x_d})=(e^{2\pi i\sum_{j=1}^dA_{1,j}x_j},\dots,e^{2\pi i\sum_{j=1}^dA_{d,j}x_j}).$$
In particular, when $d=1$ and $A=2$, the morphism $\hat\alpha_H$ on $\bt$ is $\hat\alpha_H(z)=z^2$. 

\end{example}

\begin{proposition}\label{prop6_2}
Let $H$ be a subgroup of $B$ satisfying \eqref{eqdil1}.

(i) The operators $u$ and $t_h$ defined in \eqref{eqdil4} satisfy the following relation
\begin{equation}\label{eqdil5}
ut_h=t_{\alpha(h)}u.
\end{equation}
Moreover the operators $t_h$ are unitary for all $h\in H$, and, $u$ is an isometry. 

(ii) Let 
\begin{equation}\label{eqdil6}
W:L^2(\hat H,\mu_H)\rightarrow L^2(K,\mu),\quad Wf=f\circ\theta_0,
\end{equation}
and consider the representation $U$ of $G=B\rtimes_\alpha\bz$ given in \eqref{eq39}.
Then $W$ is an isometry with the property that
\begin{equation}\label{eqdil7}
Wu=U(1,0)W,\quad Wt_h=U(0,h)W,\quad(h\in H).
\end{equation}
Moreover $U$ is the minimal unitary dilation of the isometry $u$, i.e., 
\begin{equation}\label{eqdil8}
\bigcup_{j\geq0}U(-j,0)\left[WL^2(\hat H,\mu_H)\right]\mbox{ is dense in }L^2(K,\mu).
\end{equation}

\end{proposition}

\begin{proof}
(i) Note that $t_h$ is a multiplication operator by the character $h$ in $\widehat{\hat H}=H$. Therefore it is a unitary operator. With \eqref{eqdil2}, \eqref{eqdil3}, and the invariance of $\mu$ under $\hat\alpha$ \eqref{eq34},
$$\int_{\hat H}f\circ\hat\alpha_H\,d\mu_H=\int_{\hat B}f\circ\hat\alpha_H\circ\theta_0\,d\mu=\int_{\hat B}f\circ\theta_0\circ\hat\alpha\,d\mu=\int_{\hat B}f\circ\theta_0\,d\mu=\int_{\hat H}f\,d\mu_H.$$

So $u$ is an isometry.
The relation \eqref{eqdil5} follows from a direct computation.

(ii) Equation \eqref{eqdil2} shows that $W$ is an isometry. Equations \eqref{eqdil7} follow from a direct computation (recall that $\theta_0=\hat i$).

For the \eqref{eqdil8}, we note that, if $h\in H$ is regarded as a character on $\hat H$, then
$Wh=h\circ\theta_0\in WL^2(\hat H,\mu_H)$. We have for $j\geq 0$,  $\chi$, $$(U(-j,0)(h\circ\theta_0))(\chi)=h(\theta_0\hat\alpha^{-j}(\chi))=(\hat\alpha^{-j}(\chi))(h)=(\ast),$$ since $\theta_0$ is the dual of the inclusion. Then $(\ast)=(\alpha^{-j}(h))(\chi)$. Therefore 
$U(-j,0)Wh=\alpha^{-j}h$ seen as a character on $\hat B$. 

But with \eqref{eqdil1} this means that all the characters on $\hat B$ are contained in the left-hand side of \eqref{eqdil8} so, by the Stone-Weierstrass theorem, it must be dense.
\end{proof}

\section{Wavelet questions and answers for the group $G=B\rtimes_\alpha\bz$}

    In this section we answer two natural wavelet questions which present
themselves in connection with our representations of $G = B \rtimes_\alpha\bz$. We
explore them for the class of examples in Example \ref{ex2_4}.
    Within harmonic analysis of wavelet bases it is natural to relax the
rather strict requirement that a wavelet basis be orthonormal (i.e., an
ONB). Two alternative basis possibilities are called Bessel wavelets and
frame wavelets, and we refer to the survey \cite{Chr01} for definitions and
motivation.

\begin{definition}\label{defframe}
Let $\{e_i\}_{i\in I}$ be a family of vectors in a Hilbert space $\mathcal H$. Then $\{e_i\}_{i\in I}$ is a {\it Bessel }sequence if 
there exists $M>0$ such that
$$\sum_{i\in I}|\ip{f}{e_i}|^2\leq M\|f\|^2,\quad(f\in \mathcal H).$$
The family $\{e_i\}_{i\in I}$ is called a {\it frame } if there exist constants $m,M>0$ such that
$$m\|f\|^2\leq\sum_{i\in I}|\ip{f}{e_i}|^2\leq M\|f\|^2,\quad(f\in \mathcal H).$$
\end{definition}

\begin{definition}
Consider the group in Example \ref{ex2_4}, $B=\bz_d[A^{-1}]$, and $\alpha(b)=Ab$, $b\in B$. Let $U$ be a representation of $G=B\rtimes_\alpha\bz$ on a Hilbert space $\mathcal H$. A {\it Bessel/frame/orthonormal wavelet} for the representation $U$ is a finite family $\{\psi_1,\dots,\psi_q\}$ of vectors in $\mathcal H$ such that 
$$\{U(j,\alpha^j(k))\psi_i\,|j\in\bz, k\in\bz^d, i\in\{1,\dots,q\}\}$$
is a Bessel sequence/frame/orthonormal basis for $\mathcal H$.
\end{definition}

\begin{proposition}\label{propnobesselind}
Consider the groups of Example \ref{ex2_4}, $B=\bz_d[A^{-1}]$, $\alpha(b)=Ab$ , $b\in B$. Let $\chi\in \hat B$. The representation 
$\Ind_B^G(\chi)$ has no Bessel wavelets.
\end{proposition}

\begin{proof}
We can consider the equivalent form of $\Ind_B^G(\chi)$ on $l^2(\bz)$ given in Lemma \ref{lem3_1} and Theorem \ref{th3_2}. Suppose $\{\psi_1,\dots,\psi_q\}$ is a Bessel wavelet for $\Ind_B^G(\chi)$. Then 
$$\sum_{i=1}^q\sum_{j\in\bz,k\in\bz^d}|\ip{U(j,\alpha^{j}(k))\psi_i}{\delta_0}|^2\leq M\|\delta_0\|^2=M,$$
for some constant $M>0$. But then
$$M\geq\sum_{i=1}^q\sum_{j\in\bz,k\in\bz^d}|\ip{\psi_i}{U(-j,-k)\delta_0}|^2=\sum_{i=1}^q\sum_{j\in\bz,k\in\bz^d}|\ip{\psi_i}{D_\chi(-k)T_{-j}\delta_0}|^2=$$
$$\sum_{i=1}^q\sum_{j\in\bz,k\in\bz^d}|\ip{\psi_i}{\chi(\alpha^j(k))\delta_j}|^2=\sum_{i=1}^q\sum_{k\in\bz}\sum_{j\in\bz}^d|\ip{\psi_i}{\delta_j}|^2=\sum_{i=1}^q\sum_{k\in\bz^d}\|\psi_i\|^2=\infty.$$
This contradiction proves the result.
\end{proof}

\begin{remark}\label{reml2}
For the group $B=\bz[1/2]$ in Example \ref{ex2_3} with $\alpha(x)=2x$ for $x\in B$, the representation $U$ in Theorem \ref{th4_5} ``almost'' has orthogonal wavelets. More precisely $\psi=e_1=$the character $1\in\bz[1/2]$ on $\hat B=\widehat{\bz[1/2]}$ has the property that its wavelet family $\{U(j,\alpha^{j}k)\psi\,|\,j,k\in\bz\}$ is orthonormal, and spans a space of codimension $1$ (see \cite{Dut06}).
\end{remark}

\begin{proposition}\label{propnoframe}
Consider the group in Example \ref{ex2_4}, $B=\bz_d[A^{-1}]$, $\alpha(b)=Ab$, $b\in B$. Then the right regular representation in Theorem \ref{th5_3} does not have frame wavelets.
\end{proposition}

\begin{proof} See also \cite{DaLa98} for the analysis of representations that have or do not have wavelets.

Suppose $\{\psi_1,\dots,\psi_q\}$ is a frame wavelet for the right regular representation $R$ on $l^2(G)$. 
Then 
$$m\|f\|^2\leq\sum_{i=1}^q\sum_{j\in\bz,k\in\bz^d}|\ip{R(j,\alpha^{j}(k))\psi_i}{f}|^2\leq M\|f\|^2,\quad(f\in l^2(G))$$
for some constants $m,M>0$. 
This implies that for all $g\in G$, if $f=\delta_g$ is a canonical vector,
\begin{equation}\label{eqnofr}
m\leq\sum_{i=1}^q\sum_{j\in\bz,k\in\bz^d}|\ip{\psi_i}{R(-j,-k)\delta_g}|^2=
\sum_{i=1}^q\sum_{j\in\bz,k\in\bz^d}|\ip{\psi_i}{\delta_{g(-j,-k)^{-1}}}|^2=\sum_{i=1}^q\sum_{j\in\bz,k\in\bz^d}|\ip{\psi_i}{\delta_{g(j,\alpha^j(k))}}|^2.
\end{equation}

Since $\sum_{i=1}^q\|\psi_i\|^2<\infty$, there exists a finite set $F=\{(j_1,k_1),\dots,(j_n,k_n)\}$ in $G$, such that 
\begin{equation}\label{eqnofr2}
\sum_{i=1}^q\sum_{h\in G\setminus F}|\ip{\psi_i}{\delta_h}|^2<m.
\end{equation}
Pick an element $k_0\in\bz_d[A^{-1}]$ such that $k_0+\alpha^{j_l}(k)\neq k_l$ for all $k\in\bz^d$, $l\in\{1,\dots,n\}$. This can be done by choosing $k_0=A^{-r}b_0$ with $b_0\in\bz$ and $r$ big enough (bigger than all $j_l$). 

Let $g=(0,k_0)$. Then, for $j\in\bz$ and $k\in\bz^d$, we have $h:=g(j,\alpha^j(k))=(0,k_0)(j,\alpha^j(k))=(j,k_0+\alpha^j(k))\not\in F$. 
But this, with \eqref{eqnofr2}, implies that the last term in \eqref{eqnofr} is strictly less than $m$. 

The contradiction proves the result.
\end{proof}

\section{Concluding remarks}\label{concluding}
  In our present analysis of the family of discrete groups $G = B \rtimes_\alpha \bz$  , we were motivated by certain problems in wavelets. As noted in section \ref{plancherel} above, the issue of when there is a measurable cross-section in the sense of Mackey \cite{Mac49, Mac63, Mac76} affects both the applications of this non-abelian harmonic analysis, as well as the entire structure of decompositions of data (in the form of functions) into its spectral constituents. For general references to cross-sections and operator algebras, see for example \cite{Mac63}, \cite{Arv76} and \cite{KaRi86}.

     In section \ref{plancherel}, we showed that the adjoint representation for the normal abelian subgroup B in G is simply the action of $\alpha$, and that the co-adjoint action of $G$ on $K = \hat B$ (Pontryagin dual) is the action by the dual automorphism $\hat \alpha$. Our Proposition \ref{prop3_13} states that this co-adjoint action does not have a measurable cross-section. For the benefit of readers, and for comparison, we outline below the simplest case \cite{Mac49} of a standard Mackey construction, motivated by the Stone-von Neumann uniqueness theorem, for which there {\it is} an obvious smooth cross-section.

It is interesting to compare our study of the representations of the class of discrete groups $G=B\rtimes_\alpha\bz$ with that of the Heisenberg group.
\begin{equation}\label{eq8.1}
G_{He}=\left\{\left.\left(\begin{array}{ccc}
1&a&c\\
0&1&b\\
0&0&1\end{array}\right)\right| a,b,c\in\br\right\}
\end{equation}
which was the first one Mackey considered \cite{Mac49} in connection with orbit theory.

We will do this by constrasting the Plancherel formula for $G_{He}$ with that for $G$; see Theorem \ref{th5_3} and Proposition \ref{prop3_13} (absence of measurable cross-sections) above. 

Note that $G_{He}$ is a copy of $\br^3$, and that the group law may be written in the form
\begin{equation}\label{eq8.2}
(a,b,c)(a',b',c')=(a+a',b+b',c+c'+ab').
\end{equation}
Set $B_{He}:=\{(0,b,c)\,|\,b,c\in\br\}\cong \br^2$ and note that $B_{He}$ is then a normal subgroup in $G_{He}$. If $g=(\alpha,\beta,\gamma)\in G_{He}$ is given, the adjoint action ${\operatorname*{Ad}}_g(\cdot)=g\cdot g^{-1}$ on $B_{He}$ is 
\begin{equation}\label{eq8.3}
g(0,b,c)g^{-1}=(0,b,c+\alpha b).
\end{equation}
With the standard Fourier duality  $\hat\B_{He}\cong \br^2$, we then get the Mackey co-adjoint action on $\hat B_{He}$ in the following form 
\begin{equation}\label{eq8.4}
(\xi_1,\xi_2)\stackrel{{\operatorname*{Ad}}^*_g}{\mapsto}(\xi_1+\alpha\xi_2,\xi_2)
\end{equation}
Hence there are two cases for points $\xi=(\xi_1,\xi_2)\in\hat B_{He}$.

{\it Orbit picture:} (i) If $\xi_2\neq0$, then $\mathcal O_{{\operatorname*{Ad}}_G^*}(\xi)=\br\times\{\xi_2\}$, i.e., horizontal lines.

(ii) If $\xi_2=0$, then $\mathcal O_{{\operatorname*{Ad}}_G^*}(\xi)=\{\xi\}$, i.e., the singleton $(\xi_1,0)$.

As a measurable cross-section, we may simply take a vertical line, for example $\{0\}\times\br$. It is known \cite{Mac49,Mac76} that this accounts for the Plancherel formula for $G_{He}$. Each co-adjoint orbit corresponds to a monomially induced representation $B_{He}\rightarrow G_{He}$; and the direct integral yields the regular representation of $G_{He}$.

\begin{acknowledgements}
We wish to our express our deepest thanks Professors Bob Doran, Cal Moore, and Bob Zimmer for organizing an AMS Special Session in New Orleans, January 7-8, honoring the memory of George Mackey. The paper grew out of material in a talk by the second named author at the special session, but new collaborative results between the co-authors emerged since. We have had fruitful discussions with the participants at the session; and since New Orleans with other colleagues. We want especially to thank Professors Judy Packer and Deguang Han for many very helpful discussions.
\end{acknowledgements}

\bibliographystyle{alpha}
\bibliography{dualbs}

\begin{thebibliography}{BJMP05}

\bibitem[Arv76]{Arv76}
William Arveson.
\newblock {\em An invitation to {$C\sp*$}-algebras}.
\newblock Springer-Verlag, New York, 1976.
\newblock Graduate Texts in Mathematics, No. 39.

\bibitem[BJ91]{BrJo91}
Berndt Brenken and Palle E.~T. Jorgensen.
\newblock A family of dilation crossed product algebras.
\newblock {\em J. Operator Theory}, 25(2):299--308, 1991.

\bibitem[BJMP05]{BJMP05}
Lawrence Baggett, Palle Jorgensen, Kathy Merrill, and Judith Packer.
\newblock A non-{MRA} {$C\sp r$} frame wavelet with rapid decay.
\newblock {\em Acta Appl. Math.}, 89(1-3):251--270 (2006), 2005.

\bibitem[Chr01]{Chr01}
Ole Christensen.
\newblock Frames, {R}iesz bases, and discrete {G}abor/wavelet expansions.
\newblock {\em Bull. Amer. Math. Soc. (N.S.)}, 38(3):273--291 (electronic),
  2001.

\bibitem[CMW84]{CMW84}
Ra{\'u}l~E. Curto, Paul~S. Muhly, and Dana~P. Williams.
\newblock Cross products of strongly {M}orita equivalent {$C\sp{\ast}
  $}-algebras.
\newblock {\em Proc. Amer. Math. Soc.}, 90(4):528--530, 1984.

\bibitem[Dau95]{Dau95}
Ingrid Daubechies.
\newblock Wavelets and other phase space localization methods.
\newblock In {\em Proceedings of the International Congress of Mathematicians,
  Vol.\ 1, 2 (Z\"urich, 1994)}, pages 56--74, Basel, 1995. Birkh\"auser.

\bibitem[DJ06]{DuJo06}
Dorin~Ervin Dutkay and Palle E.~T. Jorgensen.
\newblock Methods from multiscale theory and wavelets applied to nonlinear
  dynamics.
\newblock In {\em Wavelets, multiscale systems and hypercomplex analysis},
  volume 167 of {\em Oper. Theory Adv. Appl.}, pages 87--126. Birkh\"auser,
  Basel, 2006.

\bibitem[DJ07]{DuJo07}
Dorin~Ervin Dutkay and Palle E.~T. Jorgensen.
\newblock Disintegration of projective measures.
\newblock {\em Proc. Amer. Math. Soc.}, 135(1):169--179 (electronic), 2007.

\bibitem[DL98]{DaLa98}
Xingde Dai and David~R. Larson.
\newblock Wandering vectors for unitary systems and orthogonal wavelets.
\newblock {\em Mem. Amer. Math. Soc.}, 134(640):viii+68, 1998.

\bibitem[DLS98]{DLS98}
Xingde Dai, David~R. Larson, and Darrin~M. Speegle.
\newblock Wavelet sets in {$\mathbb R\sp n$}. {II}.
\newblock In {\em Wavelets, multiwavelets, and their applications (San Diego,
  CA, 1997)}, volume 216 of {\em Contemp. Math.}, pages 15--40. Amer. Math.
  Soc., Providence, RI, 1998.

\bibitem[Dut04]{Dut04}
Dorin~Ervin Dutkay.
\newblock The local trace function of shift invariant subspaces.
\newblock {\em J. Operator Theory}, 52(2):267--291, 2004.

\bibitem[Dut05]{Dut05}
Dorin~Ervin Dutkay.
\newblock Some equations relating multiwavelets and multiscaling functions.
\newblock {\em J. Funct. Anal.}, 226(1):1--20, 2005.

\bibitem[Dut06]{Dut06}
Dorin~Ervin Dutkay.
\newblock Low-pass filters and representations of the {B}aumslag {S}olitar
  group.
\newblock {\em Trans. Amer. Math. Soc.}, 358(12):5271--5291 (electronic), 2006.

\bibitem[HR63]{HeRo63}
Edwin Hewitt and Kenneth~A. Ross.
\newblock {\em Abstract harmonic analysis. {V}ol. {I}: {S}tructure of
  topological groups. {I}ntegration theory, group representations}.
\newblock Die Grundlehren der mathematischen Wissenschaften, Bd. 115. Academic
  Press Inc., Publishers, New York, 1963.

\bibitem[Jor88]{Jor88}
Palle E.~T. Jorgensen.
\newblock {\em Operators and representation theory}, volume 147 of {\em
  North-Holland Mathematics Studies}.
\newblock North-Holland Publishing Co., Amsterdam, 1988.
\newblock Canonical models for algebras of operators arising in quantum
  mechanics, Notas de Matem\'atica [Mathematical Notes], 120.

\bibitem[Jor06]{Jor06}
Palle E.~T. Jorgensen.
\newblock {\em Analysis and probability: wavelets, signals, fractals}, volume
  234 of {\em Graduate Texts in Mathematics}.
\newblock Springer, New York, 2006.

\bibitem[KR86]{KaRi86}
Richard~V. Kadison and John~R. Ringrose.
\newblock {\em Fundamentals of the theory of operator algebras. {V}ol. {II}},
  volume 100 of {\em Pure and Applied Mathematics}.
\newblock Academic Press Inc., Orlando, FL, 1986.
\newblock Advanced theory.

\bibitem[KTW85]{KTW85}
Shinz{\=o} Kawamura, Jun Tomiyama, and Yasuo Watatani.
\newblock Finite-dimensional irreducible representations of {$C\sp
  \ast$}-algebras associated with topological dynamical systems.
\newblock {\em Math. Scand.}, 56(2):241--248, 1985.

\bibitem[LPT01]{LPT01}
Lek-Heng Lim, Judith~A. Packer, and Keith~F. Taylor.
\newblock A direct integral decomposition of the wavelet representation.
\newblock {\em Proc. Amer. Math. Soc.}, 129(10):3057--3067 (electronic), 2001.

\bibitem[Mac49]{Mac49}
George~W. Mackey.
\newblock A theorem of {S}tone and von {N}eumann.
\newblock {\em Duke Math. J.}, 16:313--326, 1949.

\bibitem[Mac63]{Mac63}
George~W. Mackey.
\newblock Infinite-dimensional group representations.
\newblock {\em Bull. Amer. Math. Soc.}, 69:628--686, 1963.

\bibitem[Mac76]{Mac76}
George~W. Mackey.
\newblock {\em The theory of unitary group representations}.
\newblock University of Chicago Press, Chicago, Ill., 1976.
\newblock Based on notes by James M. G. Fell and David B. Lowdenslager of
  lectures given at the University of Chicago, Chicago, Ill., 1955, Chicago
  Lectures in Mathematics.

\bibitem[Mer05]{Mer05}
Kathy Merrill.
\newblock Constructing wavelets from generalized filters.
\newblock In {\em European women in mathematics---Marseille 2003}, volume 135
  of {\em CWI Tract}, pages 1--9. Centrum Wisk. Inform., Amsterdam, 2005.

\bibitem[MV00]{MaVa00}
Florian Martin and Alain Valette.
\newblock Markov operators on the solvable {B}aumslag-{S}olitar groups.
\newblock {\em Experiment. Math.}, 9(2):291--300, 2000.

\bibitem[{\O}rs79]{Ors79}
Bent {\O}rsted.
\newblock Induced representations and a new proof of the imprimitivity theorem.
\newblock {\em J. Funct. Anal.}, 31(3):355--359, 1979.

\bibitem[Pet83]{Pet83}
Karl Petersen.
\newblock {\em Ergodic theory}, volume~2 of {\em Cambridge Studies in Advanced
  Mathematics}.
\newblock Cambridge University Press, Cambridge, 1983.

\bibitem[Rud62]{Rud62}
Walter Rudin.
\newblock {\em Fourier analysis on groups}.
\newblock Interscience Tracts in Pure and Applied Mathematics, No. 12.
  Interscience Publishers (a division of John Wiley and Sons), New York-London,
  1962.

\bibitem[Tho64]{Tho64}
Elmar Thoma.
\newblock \"{U}ber unit\"are {D}arstellungen abz\"ahlbarer, diskreter
  {G}ruppen.
\newblock {\em Math. Ann.}, 153:111--138, 1964.

\bibitem[Wal82]{Wal82}
Peter Walters.
\newblock {\em An introduction to ergodic theory}, volume~79 of {\em Graduate
  Texts in Mathematics}.
\newblock Springer-Verlag, New York, 1982.

\bibitem[Wil82]{Wil82}
Dana~P. Williams.
\newblock Transformation group {$C\sp{\ast} $}-algebras with {H}ausdorff
  spectrum.
\newblock {\em Illinois J. Math.}, 26(2):317--321, 1982.

\end{thebibliography}

\end{document}